\newcommand{\xBc}{\langle}
\newcommand{\xBe}{\rangle}

\newcommand{\xbO}{\Omega}
\newcommand{\xbP}{\Pi}

\newcommand{\xba}{\alpha}
\newcommand{\xbb}{\beta}

\newcommand{\xbe}{\in}
\newcommand{\xbf}{\phi}
\newcommand{\xbg}{\gamma}

\newcommand{\xbk}{\kappa}

\newcommand{\xbm}{\mu}

\newcommand{\xbo}{\omega}

\newcommand{\xbq}{\psi}

\newcommand{\xbs}{\sigma}

\newcommand{\xCI}{{\Big(}}
\newcommand{\xCJ}{{\Big)}}
\newcommand{\xCK}{\times}

\newcommand{\xCN}{\neg}
\newcommand{\xCO}{ }
\newcommand{\xCQ}{\emptyset}

\newcommand{\xCf}{\hspace{0.1em}}

\newcommand{\xcA}{\forall}

\newcommand{\xcE}{\exists}

\newcommand{\xcH}{\not\Rightarrow}
\newcommand{\xcI}{\not\Leftarrow}
\newcommand{\xcJ}{\not\Leftrightarrow}

\newcommand{\xcL}{\not\vdash}

\newcommand{\xcN}{\hspace{0.2em}\not\sim\hspace{-0.9em}\mid\hspace{0.8em}}

\newcommand{\xcS}{\bigcap}
\newcommand{\xcT}{\bot}

\newcommand{\xcV}{\bigcup}

\newcommand{\xcX}{\Box}

\newcommand{\xca}{\infty}
\newcommand{\xcb}{\subset}
\newcommand{\xcc}{\subseteq}
\newcommand{\xcd}{\supseteq}
\newcommand{\xce}{\not\in}
\newcommand{\xcf}{\supset}
\newcommand{\xcg}{\geq}
\newcommand{\xch}{\Rightarrow}
\newcommand{\xci}{\Leftarrow}
\newcommand{\xcj}{\Leftrightarrow}
\newcommand{\xck}{\leq}
\newcommand{\xcl}{\vdash}
\newcommand{\xcm}{\models}
\newcommand{\xcn}{\hspace{0.2em}\sim\hspace{-0.9em}\mid\hspace{0.58em}}

\newcommand{\xco}{\vee}
\newcommand{\xcp}{\rightarrow}

\newcommand{\xcr}{\leftrightarrow}
\newcommand{\xcs}{\cap}
\newcommand{\xcu}{\wedge}
\newcommand{\xcv}{\cup}
\newcommand{\xcx}{\Diamond}

\newcommand{\xcz}{\Box}

\newcommand{\xDC}{\hspace{2em}}
\newcommand{\xDH}{\item }

\newcommand{\xDd}{\equiv}

\newcommand{\xdC}{\mbox{\boldmath$C$}}
\newcommand{\xdD}{\mbox{\boldmath$D$}}

\newcommand{\xda}{{\cal A}}

\newcommand{\xdl}{{\cal L}}
\newcommand{\xdm}{{\cal M}}

\newcommand{\xdp}{{\cal P}}

\newcommand{\xdu}{{\cal U}}

\newcommand{\xdx}{{\cal X}}
\newcommand{\xdy}{{\cal Y}}
\newcommand{\xdz}{{\cal Z}}

\newcommand{\xEH}{ & }
\newcommand{\xEI}{\begin{itemize}}
\newcommand{\xEJ}{\end{itemize}}
\newcommand{\xEP}{ \\ }

\newcommand{\xEd}{\neq}

\newcommand{\xEh}{\begin{enumerate}}

\newcommand{\xEj}{\end{enumerate}}

\newcommand{\xeb}{\prec}
\newcommand{\xec}{\preceq}

\newcommand{\xee}{\succ}

\newcommand{\xex}{\upharpoonright}

\newcommand{\xFO}{\parallel}

\newcommand{\Xl}{\ldots}

\newcommand{\ol}{\overline}

\newcommand{\xssc}{\scriptsize}

\newcommand{\bl}{\begin{lemma} \rm}
\newcommand{\el}{\end{lemma}}
\newcommand{\br}{\begin{remark} \rm}
\newcommand{\er}{\end{remark}}
\newcommand{\be}{\begin{example} \rm}
\newcommand{\ee}{\end{example}}
\newcommand{\bco}{\begin{corollary} \rm}
\newcommand{\eco}{\end{corollary}}
\newcommand{\bc}{\begin{claim} \rm}
\newcommand{\ec}{\end{claim}}
\newcommand{\bfa}{\begin{fact} \rm}
\newcommand{\efa}{\end{fact}}
\newcommand{\bp}{\begin{proposition} \rm}
\newcommand{\ep}{\end{proposition}}
\newcommand{\bd}{\begin{definition} \rm}
\newcommand{\ed}{\end{definition}}
\newcommand{\bcs}{\begin{construction} \rm}
\newcommand{\ecs}{\end{construction}}
\newcommand{\bcd}{\begin{condition} \rm}
\newcommand{\ecd}{\end{condition}}
\newcommand{\bt}{\begin{theorem} \rm}
\newcommand{\et}{\end{theorem}}
\newcommand{\bn}{\begin{notation} \rm}
\newcommand{\en}{\end{notation}}
\newcommand{\bfi}{\begin{bild} \rm}
\newcommand{\efi}{\end{bild}}
\newcommand{\bsta}{\begin{statement} \rm}
\newcommand{\esta}{\end{statement}}
\newcommand{\bcom}{\begin{comment} \rm}
\newcommand{\ecom}{\end{comment}}
\newcommand{\bdia}{\begin{diagram} \rm}
\newcommand{\edia}{\end{diagram}}

\newcommand{\bfc}{\begin{figure}[htb] \begin{center}}
\newcommand{\efc}{\end{center} \end{figure}}

\documentclass{article}

\usepackage{amssymb,latexsym,epic,eepic,rotating}

\title{
Critical analysis of
the Carmo-Jones system of
Contrary-to-Duty obligations
\thanks{
paper 358
}
}

\author{Dov M Gabbay
\thanks{
Dov.Gabbay@kcl.ac.uk, www.dcs.kcl.ac.uk/staff/dg
} \\
King's College, London
\thanks{
Department of Computer Science, King's College London, Strand,
London WC2R 2LS, UK
} \\ \\
Karl Schlechta
\thanks{
ks@cmi.univ-mrs.fr, karl.schlechta@web.de, http://www.cmi.univ-mrs.fr/ $\sim$ ks
} \\
Laboratoire d'Informatique Fondamentale de Marseille
\thanks{
UMR 6166, CNRS and Universit\'{e} de Provence,
Address: CMI, 39, rue Joliot-Curie, F-13453 Marseille Cedex 13, France
}
}

\begin{document}

\newtheorem{lemma}{Lemma}[section]
\newtheorem{theorem}[lemma]{Theorem}
\newtheorem{proposition}[lemma]{Proposition}
\newtheorem{corollary}[lemma]{Corollary}
\newtheorem{claim}[lemma]{Claim}
\newtheorem{fact}[lemma]{Fact}
\newtheorem{remark}[lemma]{Remark}
\newtheorem{definition}{Definition}[section]
\newtheorem{construction}{Construction}[section]
\newtheorem{condition}{Condition}[section]
\newtheorem{example}{Example}[section]
\newtheorem{notation}{Notation}[section]
\newtheorem{bild}{Figure}[section]
\newtheorem{comment}{Comment}[section]
\newtheorem{statement}{Statement}[section]
\newtheorem{diagram}{Diagram}[section]

\renewcommand{\labelenumi}
  {(\arabic{enumi})}
\renewcommand{\labelenumii}
  {(\arabic{enumi}.\arabic{enumii})}
\renewcommand{\labelenumiii}
  {(\arabic{enumi}.\arabic{enumii}.\arabic{enumiii})}
\renewcommand{\labelenumiv}
  {(\arabic{enumi}.\arabic{enumii}.\arabic{enumiii}.\arabic{enumiv})}

\maketitle

\setcounter{secnumdepth}{3}
\setcounter{tocdepth}{3}

\begin{abstract}

This paper offers a technical analysis of the contrary to duty system proposed
in Carmo-Jones. We offer analysis/simplification/repair of their system and
compare it with our own related system.

\end{abstract}

\tableofcontents

%
%
%
\section{
Introduction
}

\label{Section Introduction}

The present paper was inspired by the important paper  \cite{CJ02}
of J.Carmo and A.Jones on contrary-to-duties.

In that paper Carmo and Jones present a logical system designed to
solve many of the current puzzles of contrary-to-duties. They propose
a system with the unary connective $O(B)$ and the binary
connective
$O(B/A)$ and using these connectives give a detailed case analysis
of several contrary-to-duty paradoxes.

Gabbay, in his paper  \cite{Gab08} proposed a reactive Kripke
semantics
approach to contrary-to-duties and made use of the Carmo and Jones
paper to draw upon examples and analysis. Gabbay promised in his paper
an analysis of the Carmo--Jones approach and a comparison
with his own paper. Meanwhile Gabbay and Schlechta developed the reactive
and hierarchical approach to conditionals  \cite{GS08d}
as well as a general road map paper for preferential semantics  \cite{GS08c}
and armed with this new arsenal of methods (Carmo--Jones paper was
written 10 years ago), we believe we can give a preferential analysis of
the
Carmo--Jones
paper.

Our comments are strictly mathematical. Our own
philosophical approach is outlined in
 \cite{GS08g}.

\section{
The Carmo-Jones system
}

To model a contrary to duty set of sentences given in a natural language, which
avoids paradoxes, we need a logic L and a translation from natural language
into L.
The translation must be such that whenever the original natural language set is
coherent and consistent in our common sense reading of it, its natural formal
translation in L is consistent in L
(otherwise we get what is referred to as a paradox, relative to L).
One such logic L is dyadic modal logic. We have a binary modal operator $O(B/A)$
reading $B$ is obligatory relative to a given $A$, i.e., we have multiple unary
modalities $O_A$ dependent on $A$.

Thus we have

$t \xcm O(B/A)$ iff for all $s$ such that $tR(A)s$ holds we have that $s \xcm
B$.

It stands to reason that condition (5-b) below holds for $R(A)$, namely

$tR(A)s$ implies that $s$ is in $A$ (i.e., $s \xcm A$),

We note that any correct logic L needs
axioms for combining formulas of the form

$O(B/A)$ with $O(\neg B/(A \xcu C))$.

Bearing all of the above in mind, let us examine the Carmo-Jones system.

To fix our notation etc, the following is the Carmo-Jones system,
regarded formally as a logical system with axioms and semantics as
proposed by Carmo-Jones. (We take the liberty to change notation slightly,
and will sometimes call the system CJ system.)

Alphabet:

classical propositional logic, with 5 additional modal operators:

$ \xcX_{a}$ with dual $ \xcx_{a}$ - the actually necessary/possible

$ \xcX_{p}$ with dual $ \xcx_{p}$ - the potentially necessary/possible

$O(./.)$ a dyadic deontic operator

$O_{a}(.)$ monadic deontic operator: actual obligations

$O_{p}(.)$ monadic deontic operator: potenial obligations

Semantics:

A model $ \xdm = \xBc W,av,pv,ob,V \xBe $ where

(1) $W \xEd \xCQ $

(2) $V$ an assignment function

(3) $av:W \xcp \xdp (W)$ (the actually accessible worlds) such that

(3-a) $av(w) \xEd \xCQ $

(4) $pv:W \xcp \xdp (W)$ (the potentially accessible worlds) such that

(4-a) $av(w) \xcc pv(w)$

(4-b) $w \xbe pv(w)$

(5) $ob: \xdp (W) \xcp \xdp (\xdp (W))$ - the ``morally good'' sets

such that for $X,Y,Z \xcc W$

(5-a) $ \xCQ \xce ob(X)$

(5-b) if $Y \xcs X=Z \xcs X,$ then $Y \xbe ob(X) \xcj Z \xbe ob(X)$

(5-c) if $Y,Z \xbe ob(X),$ then $Y \xcs Z \xbe ob(X)$

(5-d) if $Y \xcc X \xcc Z,$ $Y \xbe ob(X),$ then $(Z-X) \xcv Y \xbe ob(Z)$

Remark: This results in a form of the Ross paradox:
Let $X:=M($ water plants),
$Y:=M($ water plants and post letter), then
$(Z-X) \xcv Y$ is the set of models where the plants are watered and the
letter
is posted (so far ok), or the plants are $ \xCf not$ watered.
So either do both, or don't water the plants - which
does not seem a good obligation.

Validity in $w$ is defined (for fixed $ \xdm)$ inductively as follows
$(M(\xbf)$ is the set of points where $ \xbf $ holds):

$w \xcm p$ $: \xcj $ $w \xbe V(p)$

the usual conditions for classical connectives

$w \xcm \xcX_{a} \xbf $ $: \xcj $ $av(w) \xcc M(\xbf)$

$w \xcm \xcX_{p} \xbf $ $: \xcj $ $pv(w) \xcc M(\xbf)$

$w \xcm O(\xbf / \xbq)$ $: \xcj $ $M(\xbf) \xcs M(\xbq) \xEd \xCQ $
and
$ \xcA X(X \xcc M(\xbq),$ $X \xcs M(\xbf) \xEd \xCQ $ $ \xch $ $M(
\xbf) \xbe ob(X))$

$w \xcm O_{a} \xbf $ $: \xcj $ $M(\xbf) \xbe ob(av(w))$ and $av(w) \xcs
M(\xCN \xbf) \xEd \xCQ $

$w \xcm O_{p} \xbf $ $: \xcj $ $M(\xbf) \xbe ob(pv(w))$ and $pv(w) \xcs
M(\xCN \xbf) \xEd \xCQ $

Axiomatics

(A) $ \xcX_{a}$ and $ \xcX_{p}$

(1) $ \xcX_{p}$ is a normal modal operator of type KT

(2) $ \xcX_{a}$ is a normal modal operator of type KD

(3) $ \xcX_{p} \xbf $ $ \xcp $ $ \xcX_{a} \xbf $

(B) Characterisation of $O(./.)$

(4) $ \xCN O(\xcT / \xbq)$

(5) $O(\xbf / \xbq) \xcu O(\xbf' / \xbq) \xcp O(\xbf \xcu \xbf' /
\xbq)$

(6) $O(\xbf / \xbq) \xcp O(\xbf / \xbf \xcu \xbq)$ $ \xCf (SA1)$

(7) If $ \xcl \xbq \xcr \xbq',$ then $ \xcl O(\xbf / \xbq) \xcr O(
\xbf / \xbq')$

(8) If $ \xcl \xbq \xcp (\xbf \xcr \xbf'),$ then $ \xcl O(\xbf / \xbq
) \xcr O(\xbf' / \xbq)$

(C) Relationship between $O(./.)$ and $ \xcX_{p}$

(9) $ \xcx_{p}O(\xbf / \xbq) \xcp \xcX_{p}O(\xbf / \xbq)$

(10) $ \xcx_{p}(\xbq \xcu \xbq' \xcu \xbf) \xcu O(\xbf / \xbq) \xcp
O(\xbf / \xbq \xcu \xbq')$ $ \xCf (SA2)$

(D) Characterization of $O_{a}$ and $O_{p}$

(11) $O_{a} \xbf \xcu O_{a} \xbq \xcp O_{a}(\xbf \xcu \xbq)$

$O_{p} \xbf \xcu O_{p} \xbq \xcp O_{p}(\xbf \xcu \xbq)$

(E) Relationships between $O_{a}$ $(O_{p})$ and $ \xcX_{a}$ $(\xcX_{p})$

(12) $ \xcX_{a} \xbf \xcp (\xCN O_{a} \xbf \xcu \xCN O_{a} \xCN \xbf)$

$ \xcX_{p} \xbf \xcp (\xCN O_{p} \xbf \xcu \xCN O_{p} \xCN \xbf)$

(13) $ \xcX_{a}(\xbf \xcr \xbq) \xcp (O_{a} \xbf \xcr O_{a} \xbq)$

$ \xcX_{p}(\xbf \xcr \xbq) \xcp (O_{p} \xbf \xcr O_{p} \xbq)$

(F) Relationships between $O(./.),$ $O_{a}$ $(O_{p})$ and $ \xcX_{a}$ $(
\xcX_{p})$

(14) $O(\xbf / \xbq)$ $ \xcu $ $ \xcX_{a} \xbq $ $ \xcu $ $ \xcx_{a}
\xbf $ $ \xcu $ $ \xcx_{a} \xCN \xbf $ $ \xcp $ $O_{a} \xbf $

$O(\xbf / \xbq)$ $ \xcu $ $ \xcX_{p} \xbq $ $ \xcu $ $ \xcx_{p} \xbf $ $
\xcu $ $ \xcx_{p} \xCN \xbf $ $ \xcp $ $O_{p} \xbf $

(15) $O(\xbf / \xbq)$ $ \xcu $ $ \xcx_{a}(\xbf \xcu \xbq)$ $ \xcu $ $
\xcx_{a}(\xbq \xcu \xCN \xbf)$ $ \xcp O_{a}(\xbq \xcp \xbf)$

$O(\xbf / \xbq)$ $ \xcu $ $ \xcx_{p}(\xbf \xcu \xbq)$ $ \xcu $ $
\xcx_{p}(\xbq \xcu \xCN \xbf)$ $ \xcp O_{p}(\xbq \xcp \xbf)$
\subsection{
General comments
}
\subsubsection{
Methodological discussion
}

We believe that Carmo and Jones important insight was that to solve
contrary-to-duty and other Deontic paradoxes we need a wider family of
operators
capable of describing a wider context surrounding the problematic
paradoxes.
We agree with this view wholeheartedly. Gabbay's papers
 \cite{Gab08} and
 \cite{Gab08a}
use reactive semantics to create such a context and the present paper will
use hierarchical modality to create essentially the same context.
See \cite{Gab08}, Example 3.1. Also note that \cite{Gab08} contains
the following text (in the current January 2010 draft
of the paper the text is on page 47):

``We can now also understand better the approach of Carmo and Jones.
Using our terminology, they were implicitly using the cut approach
by translating into a richer language with more operators, including some dyadic
ones."

It would be useful to describe the methodology we use.

Viewed formally, we have here a logical system CJ proposed by Carmo-Jones
and a
proposed semantics $ \xdm (CJ)$ for it, intended to be applied to the
contrary-to-duties application area CTD. We want to study it and compare it with
our own methodology, and technically simplify/assist/repair/support its formal
details.

We would like to provide preferential semantics for the
Carmo Jones system. How can we do it?

Let us list the methodological parameters involved.
\subsubsection{
The semantics proposed must be compatible with the intended application.
}

This means that the spirit of the semantics must correspond to the
application.

We explain by an example. Consider modal logic S4 and assume we are trying
to
apply it to the analysis of the tenses of natural language.

The phrase ``$ \xCf A$ is true from now on'' can be modelled by $ \xcX A.$

The phrase ``John is reading now'' i.e. the progressive tense can also be
modelled as $ \xcX ($ John is reading).

Both examples give rise to modal S4. However the Kripke accessibility
relation for S4 is the semantics suitable for the ``from now on''
linguistic construction, while the McKinsey-Tarski open intervals
semantics
for S4 is more suitable for the analysis of the progressive. (Sentences
$ \xCf A$ are assigned intervals $W(A)$ and $ \xcX A$ is read as the
topological interior
of $W(A).)$

Carmo-Jones indeed offer an analysis of the compatibility of their system
in Section 6 of their paper. We will examine that.
\subsubsection{
Soundness and completeness
}

We ask whether the system is sound and complete for the semantics.
(Carmo and Jones claimed only soundness.)
If not,
what axioms do we need to add to the system or what changes do we propose
to the system to obtain correspondence? We will find that CJ is not
complete for the proposed semantics.
\subsection{
Discrepancies inside the CJ system
}

A closer look at semantics and proof theory reveals a certain asymmetry in
the
treatment of unary vs. binary obligations, and elsewhere:

 \xEh

 \xDH Unary obligations are dependent on accessibility relations $ \xCf
av$ and $ \xCf pv,$
binary ones are not. As a consequence, unary obligations depend on the
world
we are in, binary ones do not.

 \xDH Unary obligations must not be trivial, i.e. the contrary must be
possible,
binary ones can be trivial.

 \xDH (And perhaps deepest) Binary obligations postulate additional
properties
of the basic choice function ob (which makes it essentially ranked), unary
obligations need only basic properties (essentially corresponding to a not
necessarily smooth preferential relation). This property is put into the
validity condition, and not into rules as one would usually expect.

 \xDH In the validity condition for $O(B/A)$ we have $X \xcc M(A)$ and $X
\xcs M(B) \xEd \xCQ,$ in
the syntactic condition (SA2) we have $ \xcx (A \xcu B \xcu C) \xcu O(C/B)
\xco O(C/A \xcu B).$ These two
coincide only if $ \xcx $ is consistency - i.e. the underlying relation is
the
trivial universal one.

 \xDH Semantic condition 5-d) gives essentially the condition for a
preferential
structure, an analogue on the syntactical side is missing - see
Example \ref{Example 1.1} (page \pageref{Example 1.1})  below, which shows that
the axioms are not
complete for the
semantics.

 \xDH We do not quite understand the derived obligation to kill and offer
a
cigarette. We think this should rather be: $O(\xCN kill),$ $O(\xCN
offer),$
$O(offer/kill).$

 \xDH $P.$ 317, violation of $O(B/A),$ a better definition seems to be:

$m$ violates $O(B/A)$ iff in $m$ holds:

$ \xcx^{-}(O(B/A)$ $ \xcu $ $ \xcx (A \xcu B))$ $ \xcu $ $A \xcu \xCN B$

$(\xcx^{-}$ is the inverse relation).

In other words: in some antecedent, $O(B/A)$ was postulated, and $ \xCf A$
and $B$ were
possible, but now (i.e. in $m)$ $A \xcu \xCN B$ holds. (We can strengthen:
$m \xcm \xcX (A \xcu \xCN B).)$

 \xDH We also think that temporal developments and intentions should
better be
coded explicitly, as implicit coding often leads to counterintuitive
results. It is not our aim to treat such aspects here.

 \xEj
\subsection{
Incompleteness of the CJ system
}

\be

$\hspace{0.01em}$


\label{Example 1.1}

Let $ \xdl $ be defined by $p,q,$ $W:=M_{ \xdl }$ be the set of its
models.

Let $m_{1} \xcm p \xcu q,$ $m_{2} \xcm p \xcu \xCN q,$ $M_{1}:=\{m_{1}\},$
$M_{2}:=\{m_{2}\}.$

We write $M(A)$ for the set of models of $ \xCf A.$

Set $ob(M_{1}):=\{M \xcc M_{ \xdl }:M_{1} \xcc M\},$ $ob(M_{2}):=\{M \xcc
M_{ \xdl }:M_{2} \xcc M\},$ $ob(M):= \xCQ $ for all other $M.$

Let $av(w):=pv(w):=W$ for all $w \xbe W,$ i.e. both are defined by $ \xCf
wRw' $ for all $w,w'.$

Thus, $O_{a}=O_{p},$ there is only one $ \xcX,$ etc., and $M \xcm_{w}
\xcX A$ iff $ \xCf A$ is a tautology.

$M \xcm_{w}OA$ will never hold, as $av(w)=W,$ and $ob(W)= \xCQ.$

$M \xcm_{w}O(B/A)$ is independent from $w,$ so we write just $M \xcm
O(B/A).$

Suppose $M \xcm O(B/A)$ holds, then $M(A) \xcs M(B) \xEd \xCQ,$ and thus
$M(B) \xbe ob(M(A)).$ So
A has to be (equivalent to) $p \xcu q$ or $p \xcu \xCN q.$ But the only
subsets of $M(A)$ are
then $ \xCQ $ and $M(A),$ and we have $O(\xbf /p \xcu q)$ iff $ \xcl p
\xcu q \xcp \xbf,$ and
$O(\xbf /p \xcu \xCN q)$ iff $ \xcl p \xcu \xCN q \xcp \xbf.$ No other
$O(A/B)$ hold.

We check the axioms (page 293-294) of  \cite{CJ02}:

1-5 are trivial.

6. is trivial, as $M \xcm O(B/A)$ implies $ \xcl A \xcp B.$

7. is trivial.

8. Let $O(A/C),$ $ \xcl C \xcp (A \xcp B),$ then $ \xcl C \xcp A,$ so $
\xcl C \xcp B,$ so $O(B/C).$

9. trivial.

10. If $O(C/B)$ and $Con(A,B,C),$ then $ \xcl B \xcp A,$ as $B$ is
complete, so $ \xcl A \xcu B \xcr B.$

11. is void.

12.-13. trivial

14. If $O(B/A),$ then $ \xCN \xcX A.$

15. If $O(B/A),$ then $ \xcl A \xcp B,$ so $ \xcx (A \xcu \xCN B)$ is
impossible.

Thus, our example satisfies the CJ axioms.

If the system were to satisfy 5-d), then $M(p)=\{m_{1},m_{2}\} \xbe
ob(M(p)),$ and
we would have $O(p/p):$

First, $M(p) \xcs M(p) \xEd \xCQ.$ We then have to consider $X=M(p),$
$M_{1},$ $M_{2}.$
But $M(p) \xbe ob(M_{1}) \xcs ob(M_{2}) \xcs ob(M(p)),$ thus $O(p/p)$
holds.

$ \xcz $
\\[3ex]
\subsection{
Simplifications of the CJ system
}

\ee

We make now some simplifications which will help us to understand the CJ system.

(1) We assume the language is finite, thus we will not have any problems
with non-definable model sets - see e.g.  \cite{GS08c} for
an illustration of what can happen otherwise.

(2) We assume that $ob(X) \xcc \xdp (X).$ This is justified by the
following fact,
which follows immediately from the system of CJ, condition 5-b):

\bfa

$\hspace{0.01em}$


\label{Fact 1.1}

If $A \xbe ob(X),$ $A \xcc X,$ $B \xcc W-$X, then $A \xcv B \xbe ob(X).$
Conversely, if $A \xbe ob(X),$ then
$A \xcs X \xbe ob(X).$

\efa

Thus, what is outside $X,$ does not matter, and we can concentrate on the
inside of $X.$ (Of course, the validity condition has then to be modified,
$M(B) \xbe ob(X)$ will be replaced by: There is $X' \xbe ob(X),$ $X' =M(B)
\xcs X.$

By 5-c), ob is closed under finite intersection, by overall finiteness,
there is thus a smallest (by $(\xcc))$ $A \xbe ob(X).$ We call this $
\xbm (X).$ Thus,
$ \xbm (X) \xcc X,$ which is condition $(\xbm \xcc).$ (If the language
is not finite,
we would have to work with the limit version. As we work with formulas
only, this would not present a fundamental problem, see
 \cite{Sch04}.)

Let $X \xcc Z,$ $Y:= \xbm (X),$ then by 5-d) $((Z-X) \xcv Y) \xbe ob(Z),$
so $ \xbm (Z) \xcc ((Z-X) \xcv Y),$ or
$ \xbm (Z) \xcs X \xcc \xbm (X),$ which is condition $(\xbm PR)$ - see
below.

We thus have that $ \xbm $ satisfies $(\xbm \xcc)$ and $(\xbm PR),$ and
we know that this
suffices for a representation by preferential structures - see
e.g.  \cite{Sch92} and Section \ref{Section Standard} (page \pageref{Section
Standard}).

Thus, the basic choice function $ \xCf ob$ is preferential for unary $O.$

Note that $(\xbm \xcc)$ $+$ $(\xbm PR)$ imply $(\xbm OR):$ $ \xbm (X
\xcv Y) \xcc \xbm (X) \xcv \xbm (Y)$ - see
 \cite{GS08c} and
Section \ref{Section Standard} (page \pageref{Section Standard}).

When we look now at the truth conditions for $O_{a}$ and $O_{i},$ we see
that
we first go to the accessible worlds - $av(w)$ or $pv(w)$ - and check
whether
$ \xbm (av(w)) \xcc M(A)$ respectively $ \xbm (pv(w)) \xcc M(A)$ (and
whether $ \xCN A$ is possible).
Thus, in preferential terms, whether $av(w) \xcn A,$ but $av(w) \xcL A.$

The case of $O(B/A)$ is a bit more complicated and is partly dissociated
from $O_{a}$ and $O_{i}.$

We said already above that $O(B/A)$ is independent from $ \xCf av$ and $
\xCf pv,$ and from $w.$

Second, and more importantly, the condition for $O(B/A)$ implies a
converse of
$(\xbm OR)$ or $(\xbm PR):$

(1) Setting $X:=M(A),$ we have $ \xbm (M(A)) \xcc M(B),$

(2) as all sets are definable, we can choose $B$ s.t. $ \xbm (M(A))=M(B),$

(3) for $X \xcc M(A),$ we have - using (2) - $ \xbm (X) \xcc \xbm (M(A))
\xcs X$ if $X \xcs \xbm (M(A)) \xEd \xCQ.$

We thus have - if $O(B/A)$ holds - together with $(\xbm PR)$ that $(\xbm
=)$ holds, i.e.

$X \xcc Y,$ $X \xcs \xbm (Y) \xEd \xCQ $ $ \xch $ $ \xbm (X)= \xbm (Y)
\xcs X.$

By 5-a) $ \xbm (X) \xEd \xCQ,$ so $(\xbm \xCQ)$ holds, too, and by
 \cite{Sch04},
see also  \cite{GS08c} and
Section \ref{Section Standard} (page \pageref{Section Standard}),
we know that such $ \xbm $ can
be represented by a ranked
smooth structure where all elements occur in one copy only.

Thus, the basic choice function $ \xCf ob$ is ranked for binary $O(B/A).$
\subsection{
Suggested modifications of the CJ system
}

 \xEh

 \xDH We assume finiteness (see above)

 \xDH We work with the smallest element of $ob(X)$ (see above)

 \xDH We use only one accessibility relation (or operation) $ \xCf a.$
This is
justified, as we are mainly interested in formal properties here.

 \xDH We make both $O$ and $O(./.)$ dependent on $ \xCf a.$ So validity of
$O(./.)$ depends
on $w,$ too.

This eliminates one discrepancy between $O$ and $O(./.).$

 \xDH We allow both $O$ and $O(./.)$ to be trivial. We could argue here
philosophically, e.g.: if you are unable to kill your grandmother, should
you then not any longer be obliged not to kill her? (No laws for jail
inmates?) But we do this rather by laziness, to simplify the basic
machinery.

This eliminates a second discrepancy.

 \xDH We take rankedness as a basic condition for $ \xCf ob,$ so it does
not depend any
more on validity of some $O(./.).$

 \xEj

We can now describe the basic ingredients of our suggested system:

 \xEh

 \xDH We take a finite ranked structure, together with - for simplicity -
one
additional relation of accessibility.

 \xDH Binary and unary obligations will be represented the same way, i.e.
the ``best'' situations will have lowest rank.

 \xDH To correspond to the usual way of speaking in deontic logic, we
translate this into a modal language, using techniques invented by
Boutelier et al.

 \xDH

 \xEj

This results in the following system.
\section{
Our proposal for a modified CJ system
}

The following is the proposed modified CJ system.
\subsection{
Our system in a preferential framework
}

Take any system for finite ranked structures.

 \xEI

 \xDH A ranked structure is defined in
Definition \ref{Definition Pref-Str} (page \pageref{Definition Pref-Str})  and
Definition \ref{Definition Rank-Rel} (page \pageref{Definition Rank-Rel}).

 \xDH Logical conditions are defined in
Definition \ref{Definition Log-Cond} (page \pageref{Definition Log-Cond}).

 \xDH Take now a characterisation, see
Proposition \ref{Proposition Pref-Representation-Without-Ref} (page
\pageref{Proposition Pref-Representation-Without-Ref}).

 \xDH For definiteness, we choose $(\xbm \xCQ),$ $(\xbm =),$ $(\xbm
\xcc).$

 \xDH We still have to add the accessibility relation $R,$ which chooses
subsets - this is trivial, as everything is definable.

 \xEJ
\subsection{
Our system in a modal framework
}

The language of obligations has usually the flavour of modal languages,
whereas the language describing preferential structures is usually
different
in decisive aspects.

If we accept that the description of obligations is suitably given by
ranked structures, then we have ready characterizations available. So
our task will be to adapt them to fit reasonably well into a modal logic
framework. We discuss this now.

We will suppose that we have an entry point $u$ into the structure, from
which
all models are visible through relation $R,$ with modal operators $ \xcX $
and $ \xcx.$
$R$ is supposed to be transitive.

The first hurdle is to express minimality in modal terms. Boutilier and
Lamarre have shown how to do it, see
 \cite{Bou90a} and  \cite{Lam91}. (It was criticized
in  \cite{Mak93}, but this criticism does not concern our approach
as we use
different relations for accessibility and minimization.)

We introduce a new modal operator working with the minimality relation,
say we call the (irreflexive) relation $R',$ and the corresponding
operators $ \xcX' $
and $ \xcx'.$
Being a minimal model of $ \xba $ can now be expressed by $m \xcm \xba
\xcu \xCN \xcx' \xba.$

So $ \xba \xcn \xbb $ reads: $u \xcm \xcX ((\xba \xcu \xCN \xcx' \xba)
\xcp \xbb)$ - everywhere, if $m$ is a minimal
model of $ \xba,$ then $ \xbb $ holds.

$ \xCf (RatM)$ e.g. is translated to

$u$ $ \xcm $ $ \xcX \xCI (\xbf \xcu \xCN \xcx' \xbf) \xcp \xbq \xCJ $ $
\xcu $ $ \xcx \xCI \xbf \xcu \xCN \xcx' \xbf \xcu \xbq' \xCJ $ $ \xcp $
$ \xcX \xCI (\xbf \xcu \xbq' \xcu \xCN \xcx' (\xbf \xcu \xbq')) \xcp
\xbq \xCJ.$

The second hurdle is to handle subsets defined by accessibility from a
given model $m.$ In above example, all was done from $u,$ with formulas.
But we also have to make sure that we can handle expressions like
``in all best models among those accessible from $m$ $ \xbf $ holds''. The
set of
all those accessible models corresponds to some $ \xbf_{m},$ and then we
have
to choose the best among them. In particular, we have to make sure that
the axioms of our system hold not only for the models of some formulas
seen
from $u,$ but also when those formulas are defined by the set of models
accessible from some model $m.$

Let $R(m):=\{n:mRn\},$ and $ \xbm (X)$ be the minimal models of $X.$

Suppose we want to say now: If $ \xCf mRm' $ (so $R(m') \xcc R(m)$ by
transitivity), and
$R(m') \xcs \xbm (R(m)) \xEd \xCQ,$ then
$ \xbm (R(m'))=R(m') \xcs \xbm (R(m)).$ How can we express this with
modal formulas?
If we write $m \xcm \xcX \xbf,$ then we know that $ \xbf $ holds
everywhere in $R(m),$ but $ \xbf $
might not be precise enough to describe $R(m),$ e.g. $ \xbf $ might be
TRUE.

We introduce an auxiliary modal relation $R_{-}$ with operators $
\xcX_{-}$ and $ \xcx_{-}$ s.t.
$mR_{-}m' $ iff $not(mRm').$ (If $R$ is not reflexive, $R_{-}$ will not
be either, and we
change the definition accordingly. - Our notation differs from the one
of Boutilier, we chose it as we do not know how to create his symbols.)

We can now characterize $R(m)$ by $ \xbf_{m}:$ $m \xcm \xcX \xbf_{m} \xcu
\xCN \xcx_{-} \xbf_{m}$ - everywhere
$ \xbf_{m}$ holds, and at no point we cannot reach from $m,$ $ \xbf_{m}$
holds.
We can now express that $ \xbf $ holds in the minimal models of $R(m)$ by

$m \xcm (\xcX \xbf_{m} \xcu \xCN \xcx_{-} \xbf_{m}) \xcu \xcX ((\xbf_{m}
\xcu \xCN \xcx' \xbf_{m}) \xcp \xbf).$

Finally, we can express e.g. $ \xCf (AND)$

$ \xba \xcn \xbf,$ $ \xba \xcn \xbf' $ $ \xch $ $ \xba \xcn \xbf \xcu
\xbf' $

in the case where $ \xba $ is defined by some $R(m)$ as follows:

$u$ $ \xcm $ $ \xcX $ $ \xCI $ $(\xcX \xbf_{m} \xcu \xCN \xcx_{-}
\xbf_{m})$ $ \xcu $ $ \xcX ((\xbf_{m} \xcu \xCN \xcx' \xbf_{m}) \xcp
\xbf)$ $ \xcu $
$ \xcX ((\xbf_{m} \xcu \xCN \xcx' \xbf_{m}) \xcp \xbf')$ $ \xcp $

$ \xDC $ $ \xcX ((\xbf_{m} \xcu \xCN \xcx' \xbf_{m}) \xcp \xbf \xcu \xbf
')$ $ \xCJ.$
\section{
Comparison to other systems
}

We point out here the main points of  \cite{CJ02},  \cite{GS08d}, and the
present
article, which differentiate them from the others.

 \xEI

 \xDH The Carmo-Jones article

 \xEh

 \xDH It contains much material on motivation, and discussion of examples
and paradoxa.

 \xDH It gives an account of the differences between describing situations
and valid obligations.

 \xDH It presents a descriptive semantics.

 \xDH It puts the operators in the object language and uses a modal
logic language, as usual in the field.

 \xEj

 \xDH The article on $ \xda -$ranked semantics,  \cite{GS08d}:

 \xEh

 \xDH It contains a relatively exhaustive semantics for the ideal cases
in contrary-to-duty obligations.

 \xEI

 \xDH The $ \xda -$ranked semantics allows us to express that a whole
hierarchy
of obligations (if  \Xl. possible, then  \Xl.; if not, but  \Xl., then
 \Xl.; \Xl.) is satisfied, i.e. the agent ``does his best''. This
hierarchy is directly built into the semantics, which is a
multi-layered, semi-ranked structure, which can also be re-used
in other contexts.

 \xDH The article contains a sound and complete characterization of
the semantics with full proofs.

 \xDH The language is that of usual nonmonotonic logics, i.e. rules are
given in the meta-language.

 \xEJ

 \xDH Paradoxa like the Ross paradox are not treated at all, we only
treat the ideal case, and not individual obligations.

 \xDH The additional accessibility relation is added without changing the
overall language to a modal flavour.

 \xEj

 \xDH The article on the semantics of obligations,  \cite{GS08g}:

 \xEh

 \xDH In this article, we present a discussion of elementary properties
a notion of derivation of obligations should have.

 \xDH There, we are not at all concerned about more complicated
situations, involving accessibility etc.

 \xDH We also see rankedness somewhat sceptically there.

 \xEj

 \xDH The present article

 \xEh

 \xDH We work with a ranked structure describing ideal situations as
usual.

 \xDH We fully integrate the underlying logic for the ideal cases in a
modal framework, using an idea by Boutelier and Lamarre, and
extending it with a complementary relation to precisely characterize
the successor sets.

 \xEj

 \xEJ
\section{
Definitions and proofs
}

\label{Section Standard}

$ \xCO $
\index{Definition Algebraic basics}

\bd

$\hspace{0.01em}$


\label{Definition Alg-Base}

 \xEh

 \xDH

We use $ \xdp $ \index{$ \xdp $}  to denote the power set operator,
 \index{$ \xbP $}
$ \xbP \{X_{i}:i \xbe I\}$ $:=$ $\{g:$ $g:I \xcp \xcV \{X_{i}:i \xbe I\},$
$ \xcA i \xbe I.g(i) \xbe X_{i}\}$ is the general cartesian
product, $card(X)$ shall denote the cardinality of $X,$
 \index{card}
and $V$ \index{$V$}  the
set-theoretic
universe we work in - the class of all sets. Given a set of pairs $ \xdx
,$ and a
set $X,$ we denote by $ \xdx \xex X:=\{ \xBc x,i \xBe  \xbe \xdx:x \xbe X\}.$
 \index{$ \xex $}
When the context is clear, we
will sometime simply write $X$ for $ \xdx \xex X.$ (The intended use is
for preferential
structures, where $x$ will be a point (intention: a classical
propositional
model), and $i$ an index, permitting copies of logically identical
points.)

 \xDH

$A \xcc B$ will denote that $ \xCf A$ is a subset of $B$ or equal to $B,$
and $A \xcb B$ that $ \xCf A$ is
a proper subset of $B,$ likewise for $A \xcd B$ and $A \xcf B.$
 \index{$ \xcc $}
 \index{$ \xcb $}
 \index{$ \xcd $}
 \index{$ \xcf $}

Given some fixed set $U$ we work in, and $X \xcc U,$ then $ \xdC (X):=U-X$
.
 \index{$ \xdC $}

 \xDH

If $ \xdy \xcc \xdp (X)$ for some
$X,$ we say that $ \xdy $ satisfies

$(\xcs)$ \index{$(\xcs)$}  iff it is closed under finite intersections,

$(\xcS)$ \index{$(\xcS)$}  iff it is closed under arbitrary
intersections,

$(\xcv)$ \index{$(\xcv)$}  iff it is closed under finite unions,

$(\xcV)$ \index{$(\xcV)$}  iff it is closed under arbitrary unions,

$(\xdC)$ \index{$(\xdC)$}  iff it is closed under complementation,

$ \xCf (-)$ \index{$ \xCf (-)$}  iff it is closed under set difference.

 \xDH

We will sometimes write $A=B \xFO C$ for: $A=B,$ or $A=C,$ or $A=B \xcv
C.$
 \index{$ \xFO $}

 \xEj

We make ample and tacit use of the Axiom of Choice.

\ed

$ \xCO $

$ \xCO $
\index{Definition Logic, basics}

\bd

$\hspace{0.01em}$


\label{Definition Log-Base}

 \xEh

 \xDH

We work here in a classical propositional language $ \xdl,$ a theory
\index{theory}  $T$ will be
an
arbitrary set of formulas. Formulas will often be named $ \xbf,$ $ \xbq
,$ etc., theories
$T,$ $S,$ etc.

$v(\xdl)$ \index{$v(\xdl)$}  will be the set of propositional variables
of $ \xdl.$

$F(\xdl)$ \index{$F(\xdl)$}  will be the set of formulas of $ \xdl.$

$M_{ \xdl }$ \index{$M_{ \xdl }$}  will be the set of (classical) models for $
\xdl,$
$M(T)$ \index{$M(T)$}  or
$M_{T}$ \index{$M_{T}$}
is the set of models of $T,$ likewise $M(\xbf)$ \index{$M(\xbf)$}  for a
formula $ \xbf.$

 \xDH

$ \xdD_{ \xdl }$ \index{$ \xdD_{ \xdl }$}  $:=\{M(T):$ $T$ a theory in $ \xdl
\},$
the set of
$ \xCf definable$ \index{$ \xCf definable$}  model sets.

Note that, in classical propositional logic, $ \xCQ,M_{ \xdl } \xbe
\xdD_{ \xdl },$ $ \xdD_{ \xdl }$ contains
singletons, is closed under arbitrary intersections and finite unions.

An operation $f: \xdy \xcp \xdp (M_{ \xdl })$ for $ \xdy \xcc \xdp (M_{
\xdl })$ is called
$ \xCf definability$ $ \xCf preserving$ \index{$ \xCf definability$ $ \xCf
preserving$},
$ \xCf (dp)$ \index{$ \xCf (dp)$}  or
$(\xbm dp)$ \index{$(\xbm dp)$}  in short, iff for all $X \xbe \xdD_{ \xdl }
\xcs \xdy $ $f(X) \xbe \xdD_{ \xdl }.$

We will also use $(\xbm dp)$ for binary functions $f: \xdy \xCK \xdy \xcp
\xdp (M_{ \xdl })$ - as needed
for theory revision - with the obvious meaning.

 \xDH

$ \xcl $ \index{$ \xcl $}  will be classical derivability, and

$ \ol{T}:=\{ \xbf:T \xcl \xbf \},$ the closure of $T$ under $ \xcl.$

 \index{$ \ol{T}$}

 \xDH

$Con(.)$ \index{$Con(.)$}  will stand for classical consistency, so $Con(
\xbf)$ will mean that
$ \xbf $ is classical consistent, likewise for $Con(T).$ $Con(T,T')$ will
stand for
$Con(T \xcv T'),$ etc.

 \xDH

Given a consequence relation $ \xcn,$ we define

$ \ol{ \ol{T} }:=\{ \xbf:T \xcn \xbf \}.$

 \index{$ \ol{ \ol{T} }$}

(There is no fear of confusion with $ \ol{T},$ as it just is not useful to
close
twice under classical logic.)

 \xDH

$T \xco T' $ \index{$T \xco T' $}  $:=\{ \xbf \xco \xbf': \xbf \xbe T, \xbf
' \xbe T' \}.$

 \xDH

If $X \xcc M_{ \xdl },$ then $Th(X)$ \index{$Th(X)$}  $:=\{ \xbf:X \xcm
\xbf \},$ likewise for
$Th(m)$ \index{$Th(m)$}, $m \xbe M_{ \xdl }.$
($ \xcm $ \index{$ \xcm $}  will usually be classical validity.)

 \xEj

\ed

$ \xCO $

$ \xCO $
\index{Definition Logical rules}

\bd

$\hspace{0.01em}$


\label{Definition Log-Cond}

We introduce here formally a list of properties of set functions on the
algebraic side, and their corresponding logical rules on the other side.
Putting them in parallel facilitates orientation, especially when
considering representation problems.

We show, wherever adequate, in parallel the formula version
in the left column, the theory version
in the middle column, and the semantical or algebraic
counterpart in the
right column. The algebraic counterpart gives conditions for a
function $f:\xdy\xcp\xdp (U)$, where $U$ is some set, and
$\xdy\xcc\xdp (U)$.

The development in two directions, vertically with often increasing strength,
horizontally connecting proof theory with semantics motivates the presentation
in a table. The table is split in two, as one table would be too big to print.
The first table contains the basic rules, the second one those about
cumulativity and rationality.

\ed

Precise connections between the columns are given in
Proposition \ref{Proposition Alg-Log} (page \pageref{Proposition Alg-Log}).

When the formula version is not commonly used, we omit it,
as we normally work only with the theory version.

$A$ and $B$ in the right hand side column stand for
$M(\xbf)$ for some formula $\xbf$, whereas $X$, $Y$ stand for
$M(T)$ for some theory $T$.

\begin{table}[h]
\caption{Basic logical and semantic laws}
{\xssc
\begin{tabular}{|c|c|c|}

\hline

\multicolumn{3}{|c|}{Basics} \xEP

\hline

$(AND)$
\xEH
$(AND)$
\xEH
Closure under
\xEP

$ \xbf \xcn \xbq,  \xbf \xcn \xbq'   \xch $
\xEH
$ T \xcn \xbq, T \xcn \xbq'   \xch $
\xEH
finite
\xEP

$ \xbf \xcn \xbq \xcu \xbq' $
\xEH
$ T \xcn \xbq \xcu \xbq' $
\xEH
intersection
\xEP

\hline

$(OR)$ \xEH $(OR)$ \xEH $(\xbm OR)$ \xEP

$ \xbf \xcn \xbq,  \xbf' \xcn \xbq   \xch $ \xEH
$ \ol{\ol{T}} \xcs \ol{\ol{T'}} \xcc \ol{\ol{T \xco T'}} $ \xEH
$f(X \xcv Y) \xcc f(X) \xcv f(Y)$
\xEP

$ \xbf \xco \xbf' \xcn \xbq $ \xEH
\xEH
\xEP

\hline

$(wOR)$
\xEH
$(wOR)$
\xEH
$(\xbm wOR)$
\xEP

$ \xbf \xcn \xbq,$ $ \xbf' \xcl \xbq $ $ \xch $
\xEH
$ \ol{ \ol{T} } \xcs \ol{T' }$ $ \xcc $ $ \ol{ \ol{T \xco T' } }$
\xEH
$f(X \xcv Y) \xcc f(X) \xcv Y$
\xEP

$ \xbf \xco \xbf' \xcn \xbq $
\xEH
\xEH
\xEP

\hline

$(disjOR)$
\xEH
$(disjOR)$
\xEH
$(\xbm disjOR)$
\xEP

$ \xbf \xcl \xCN \xbf',$ $ \xbf \xcn \xbq,$
\xEH
$\xCN Con(T \xcv T') \xch$
\xEH
$X \xcs Y= \xCQ $ $ \xch $
\xEP

$ \xbf' \xcn \xbq $ $ \xch $ $ \xbf \xco \xbf' \xcn \xbq $
\xEH
$ \ol{ \ol{T} } \xcs \ol{ \ol{T' } } \xcc \ol{ \ol{T \xco T' } }$
\xEH
$f(X \xcv Y) \xcc f(X) \xcv f(Y)$
\xEP

\hline

$(LLE)$
\xEH
$(LLE)$
\xEH
\xEP

Left Logical Equivalence
\xEH
\xEH
\xEP

$ \xcl \xbf \xcr \xbf',  \xbf \xcn \xbq   \xch $
\xEH
$ \ol{T}= \ol{T' }  \xch   \ol{\ol{T}} = \ol{\ol{T'}}$
\xEH
trivially true
\xEP

$ \xbf' \xcn \xbq $ \xEH \xEH \xEP

\hline

$(RW)$ Right Weakening
\xEH
$(RW)$
\xEH
upward closure
\xEP

$ \xbf \xcn \xbq,  \xcl \xbq \xcp \xbq'   \xch $
\xEH
$ T \xcn \xbq,  \xcl \xbq \xcp \xbq'   \xch $
\xEH
\xEP

$ \xbf \xcn \xbq' $
\xEH
$T \xcn \xbq' $
\xEH
\xEP

\hline

$(CCL)$ Classical Closure \xEH $(CCL)$ \xEH \xEP

\xEH
$ \ol{ \ol{T} }$ is classically
\xEH
trivially true
\xEP

\xEH closed \xEH \xEP

\hline

$(SC)$ Supraclassicality \xEH $(SC)$ \xEH $(\xbm \xcc)$ \xEP

$ \xbf \xcl \xbq $ $ \xch $ $ \xbf \xcn \xbq $ \xEH $ \ol{T} \xcc \ol{
\ol{T} }$ \xEH $f(X) \xcc X$ \xEP

\cline{1-1}

$(REF)$ Reflexivity \xEH \xEH \xEP
$ T \xcv \{\xba\} \xcn \xba $ \xEH \xEH \xEP

\hline

$(CP)$ \xEH $(CP)$ \xEH $(\xbm \xCQ)$ \xEP

Consistency Preservation \xEH \xEH \xEP

$ \xbf \xcn \xcT $ $ \xch $ $ \xbf \xcl \xcT $ \xEH $T \xcn \xcT $ $ \xch
$ $T \xcl \xcT $ \xEH $f(X)= \xCQ $ $ \xch $ $X= \xCQ $ \xEP

\hline

\xEH
\xEH $(\xbm \xCQ fin)$
\xEP

\xEH
\xEH $X \xEd \xCQ $ $ \xch $ $f(X) \xEd \xCQ $
\xEP

\xEH
\xEH for finite $X$
\xEP

\hline

\xEH $(PR)$ \xEH $(\xbm PR)$ \xEP

$ \ol{ \ol{ \xbf \xcu \xbf' } }$ $ \xcc $ $ \ol{ \ol{ \ol{ \xbf } } \xcv
\{ \xbf' \}}$ \xEH
$ \ol{ \ol{T \xcv T' } }$ $ \xcc $ $ \ol{ \ol{ \ol{T} } \xcv T' }$ \xEH
$X \xcc Y$ $ \xch $
\xEP

\xEH \xEH $f(Y) \xcs X \xcc f(X)$
\xEP

\cline{3-3}

\xEH
\xEH
$(\xbm PR')$
\xEP

\xEH
\xEH
$f(X) \xcs Y \xcc f(X \xcs Y)$
\xEP

\hline

$(CUT)$ \xEH $(CUT)$ \xEH $ (\xbm CUT) $ \xEP
$ T  \xcn \xba; T \xcv \{ \xba\} \xcn \xbb \xch $ \xEH
$T \xcc \ol{T' } \xcc \ol{ \ol{T} }  \xch $ \xEH
$f(X) \xcc Y \xcc X  \xch $ \xEP
$ T  \xcn \xbb $ \xEH
$ \ol{ \ol{T'} } \xcc \ol{ \ol{T} }$ \xEH
$f(X) \xcc f(Y)$
\xEP

\hline

\end{tabular}
}
\end{table}

\begin{table}[h]
\caption{Cumulativity and Rationality}
\tabcolsep=0.5pt
{\xssc
\begin{tabular}{|c|c|c|}

\hline

\multicolumn{3}{|c|}{Cumulativity} \xEP

\hline

$(CM)$ Cautious Monotony \xEH $(CM)$ \xEH $ (\xbm CM) $ \xEP

$ \xbf \xcn \xbq,  \xbf \xcn \xbq'   \xch $ \xEH
$T \xcc \ol{T' } \xcc \ol{ \ol{T} }  \xch $ \xEH
$f(X) \xcc Y \xcc X  \xch $
\xEP

$ \xbf \xcu \xbq \xcn \xbq' $ \xEH
$ \ol{ \ol{T} } \xcc \ol{ \ol{T' } }$ \xEH
$f(Y) \xcc f(X)$
\xEP

\cline{1-1}

\cline{3-3}

or $(ResM)$
\xEH
\xEH
$(\xbm ResM)$
\xEP

Restricted Monotony
\xEH
\xEH
$ f(X) \xcc A \xcs B \xch $
\xEP

$ T  \xcn \xba, \xbb \xch T \xcv \{\xba\} \xcn \xbb $
\xEH
\xEH
$ f(X \xcs A) \xcc B $
\xEP

\hline

$(CUM)$ Cumulativity \xEH $(CUM)$ \xEH $(\xbm CUM)$ \xEP

$ \xbf \xcn \xbq   \xch $ \xEH
$T \xcc \ol{T' } \xcc \ol{ \ol{T} }  \xch $ \xEH
$f(X) \xcc Y \xcc X  \xch $
\xEP

$(\xbf \xcn \xbq'   \xcj   \xbf \xcu \xbq \xcn \xbq')$ \xEH
$ \ol{ \ol{T} }= \ol{ \ol{T' } }$ \xEH
$f(Y)=f(X)$ \xEP

\hline

\xEH
$ (\xcc \xcd) $
\xEH
$ (\xbm \xcc \xcd) $
\xEP
\xEH
$T \xcc \ol{\ol{T'}}, T' \xcc \ol{\ol{T}} \xch $
\xEH
$ f(X) \xcc Y, f(Y) \xcc X \xch $
\xEP
\xEH
$ \ol{\ol{T'}} = \ol{\ol{T}}$
\xEH
$ f(X)=f(Y) $
\xEP

\hline

\multicolumn{3}{|c|}{Rationality} \xEP

\hline

$(RatM)$ Rational Monotony \xEH $(RatM)$ \xEH $(\xbm RatM)$ \xEP

$ \xbf \xcn \xbq,  \xbf \xcN \xCN \xbq'   \xch $ \xEH
$Con(T \xcv \ol{\ol{T'}})$, $T \xcl T'$ $ \xch $ \xEH
$X \xcc Y, X \xcs f(Y) \xEd \xCQ   \xch $
\xEP

$ \xbf \xcu \xbq' \xcn \xbq $ \xEH
$ \ol{\ol{T}} \xcd \ol{\ol{\ol{T'}} \xcv T} $ \xEH
$f(X) \xcc f(Y) \xcs X$ \xEP

\hline

\xEH $(RatM=)$ \xEH $(\xbm =)$ \xEP

\xEH
$Con(T \xcv \ol{\ol{T'}})$, $T \xcl T'$ $ \xch $ \xEH
$X \xcc Y, X \xcs f(Y) \xEd \xCQ   \xch $
\xEP

\xEH
$ \ol{\ol{T}} = \ol{\ol{\ol{T'}} \xcv T} $ \xEH
$f(X) = f(Y) \xcs X$ \xEP

\hline

\xEH
$(Log=')$
\xEH $(\xbm =')$
\xEP

\xEH
$Con(\ol{ \ol{T' } } \xcv T)$ $ \xch $
\xEH $f(Y) \xcs X \xEd \xCQ $ $ \xch $
\xEP

\xEH
$ \ol{ \ol{T \xcv T' } }= \ol{ \ol{ \ol{T' } } \xcv T}$
\xEH $f(Y \xcs X)=f(Y) \xcs X$
\xEP

\hline

\xEH
$(Log \xFO)$
\xEH $(\xbm \xFO)$
\xEP

\xEH
$ \ol{ \ol{T \xco T' } }$ is one of
\xEH $f(X \xcv Y)$ is one of
\xEP

\xEH
$\ol{\ol{T}},$ or $\ol{\ol{T'}},$ or $\ol{\ol{T}} \xcs \ol{\ol{T'}}$
\xEH
$f(X),$ $f(Y)$ or $f(X) \xcv f(Y)$
\xEP

\xEH
(by (CCL))
\xEH
\xEP

\hline

\xEH
$(Log \xcv)$
\xEH $(\xbm \xcv)$
\xEP

\xEH
$Con(\ol{ \ol{T' } } \xcv T),$ $ \xCN Con(\ol{ \ol{T' } }
\xcv \ol{ \ol{T} })$
\xEH $f(Y) \xcs (X-f(X)) \xEd \xCQ $ $ \xch $
\xEP

\xEH
$ \xch $ $ \xCN Con(\ol{ \ol{T \xco T' } } \xcv T')$
\xEH $f(X \xcv Y) \xcs Y= \xCQ$
\xEP

\hline

\xEH
$(Log \xcv')$
\xEH $(\xbm \xcv')$
\xEP

\xEH
$Con(\ol{ \ol{T' } } \xcv T),$ $ \xCN Con(\ol{ \ol{T' }
} \xcv \ol{ \ol{T} })$
\xEH $f(Y) \xcs (X-f(X)) \xEd \xCQ $ $ \xch $
\xEP

\xEH
$ \xch $ $ \ol{ \ol{T \xco T' } }= \ol{ \ol{T} }$
\xEH $f(X \xcv Y)=f(X)$
\xEP

\hline

\xEH
\xEH $(\xbm \xbe)$
\xEP

\xEH
\xEH $a \xbe X-f(X)$ $ \xch $
\xEP

\xEH
\xEH $ \xcE b \xbe X.a \xce f(\{a,b\})$
\xEP

\hline

\end{tabular}
}
\end{table}

\begin{itemize}

\item

$(PR)$ is also called $infinite$ $conditionalization$ We choose this name for
its central role for preferential structures $(PR)$ or $(\xbm PR).$

\item

The system of rules $(AND)$ $(OR)$ $(LLE)$ $(RW)$ $(SC)$ $(CP)$ $(CM)$ $(CUM)$
is also called system $P$ (for preferential). Adding $(RatM)$ gives the system
$R$ (for rationality or rankedness).

Roughly: Smooth preferential structures generate logics satisfying system
$P$, while ranked structures generate logics satisfying system $R$.

\item

A logic satisfying $(REF)$, $(ResM)$, and $(CUT)$ is called a $consequence$
$relation$.

\item

$(LLE)$ and$(CCL)$ will hold automatically, whenever we work with model sets.

\item

$(AND)$ is obviously closely related to filters, and corresponds to closure
under finite intersections. $(RW)$ corresponds to upward closure of filters.

More precisely, validity of both depend on the definition, and the
direction we consider.

Given $f$ and $(\xbm \xcc)$, $f(X)\xcc X$ generates a principal filter:
$\{X'\xcc X:f(X)\xcc X'\}$, with
the definition: If $X=M(T)$, then $T\xcn \xbf$  iff $f(X)\xcc M(\xbf)$.
Validity of $(AND)$ and
$(RW)$ are then trivial.

Conversely, we can define for $X=M(T)$

$\xdx:=\{X'\xcc X: \xcE \xbf (X'=X\xcs M(\xbf)$ and $T\xcn \xbf)\}$.

$(AND)$ then makes $\xdx$  closed under
finite intersections, and $(RW)$ makes $\xdx$  upward
closed. This is in the infinite case usually not yet a filter, as not all
subsets of $X$ need to be definable this way.
In this case, we complete $\xdx$  by
adding all $X''$ such that there is $X'\xcc X''\xcc X$, $X'\xbe\xdx$.

Alternatively, we can define

$\xdx:=\{X'\xcc X: \xcS\{X \xcs M(\xbf): T\xcn \xbf \} \xcc X' \}$.

\item

$(SC)$ corresponds to the choice of a subset.

\item

$(CP)$ is somewhat delicate, as it presupposes that the chosen model set is
non-empty. This might fail in the presence of ever better choices, without
ideal ones; the problem is addressed by the limit versions.

\item

$(PR)$ is an infinitary version of one half of the deduction theorem: Let $T$
stand for $\xbf$, $T'$ for $\xbq$, and $\xbf \xcu \xbq \xcn \xbs$,
so $\xbf \xcn \xbq \xcp \xbs$, but $(\xbq \xcp \xbs)\xcu \xbq \xcl \xbs$.

\item

$(CUM)$ (whose more interesting half in our context is $(CM)$) may best be seen
as normal use of lemmas: We have worked hard and found some lemmas. Now
we can take a rest, and come back again with our new lemmas. Adding them to the
axioms will neither add new theorems, nor prevent old ones to hold.
(This is, of course, a meta-level argument concerning an object level rule.
But also object level rules should - at least generally - have an intuitive
justification, which will then come from a meta-level argument.)

\end{itemize}

$ \xCO $

$ \xCO $

\bfa

$\hspace{0.01em}$


\label{Fact Mu-Base}

The following table is to be read as follows: If the left hand side holds
for
some
function $f: \xdy \xcp \xdp (U),$ and the auxiliary properties noted in
the middle also
hold for $f$ or $ \xdy,$ then the right hand side will hold, too - and
conversely.

``sing.'' will stand for:
``$ \xdy $ contains singletons''

\begin{table}
\caption{Interdependencies of algebraic rules}
\tabcolsep=0.5pt.
{\xssc
\begin{tabular}{|c|c|c|c|}

\hline

\multicolumn{4}{|c|}{Basics} \xEP

\hline

(1.1)
\xEH
$(\xbm PR)$
\xEH
$\xch$ $(\xcs)+(\xbm \xcc)$
\xEH
$(\xbm PR')$
\xEP

\cline{1-1}

\cline{3-3}

(1.2)
\xEH
\xEH
$\xci$
\xEH
\xEP

\hline

(2.1)
\xEH
$(\xbm PR)$
\xEH
$\xch$ $(\xbm \xcc)$
\xEH
$(\xbm OR)$
\xEP

\cline{1-1}

\cline{3-3}

(2.2)
\xEH
\xEH
$\xci$ $(\xbm \xcc)$ + $(-)$
\xEH
\xEP

\cline{1-1}

\cline{3-4}

(2.3)
\xEH
\xEH
$\xch$ $(\xbm \xcc)$
\xEH
$(\xbm wOR)$
\xEP

\cline{1-1}

\cline{3-3}

(2.4)
\xEH
\xEH
$\xci$ $(\xbm \xcc)$ + $(-)$
\xEH
\xEP

\hline

(3)
\xEH
$(\xbm PR)$
\xEH
$\xch$
\xEH
$(\xbm CUT)$
\xEP

\hline

(4)
\xEH
$(\xbm \xcc)+(\xbm \xcc \xcd)+(\xbm CUM)$
\xEH
$\xcH$
\xEH
$(\xbm PR)$
\xEP

\xEH
$+(\xbm RatM)+(\xcs)$
\xEH
\xEH
\xEP

\hline

\multicolumn{4}{|c|}{Cumulativity} \xEP

\hline

(5.1)
\xEH
$(\xbm CM)$
\xEH
$\xch$ $(\xcs)+(\xbm \xcc)$
\xEH
$(\xbm ResM)$
\xEP

\cline{1-1}

\cline{3-3}

(5.2)
\xEH
\xEH
$\xci$ (infin.)
\xEH
\xEP

\hline

(6)
\xEH
$(\xbm CM)+(\xbm CUT)$
\xEH
$\xcj$
\xEH
$(\xbm CUM)$
\xEP

\hline

(7)
\xEH
$(\xbm \xcc)+(\xbm \xcc \xcd)$
\xEH
$\xch$
\xEH
$(\xbm CUM)$
\xEP

\hline

(8)
\xEH
$(\xbm \xcc)+(\xbm CUM)+(\xcs)$
\xEH
$\xch$
\xEH
$(\xbm \xcc \xcd)$
\xEP

\hline

(9)
\xEH
$(\xbm \xcc)+(\xbm CUM)$
\xEH
$\xcH$
\xEH
$(\xbm \xcc \xcd)$
\xEP

\hline

\multicolumn{4}{|c|}{Rationality} \xEP

\hline

(10)
\xEH
$(\xbm RatM)+(\xbm PR)$
\xEH
$\xch$
\xEH
$(\xbm =)$
\xEP

\hline

(11)
\xEH
$(\xbm =)$
\xEH
$ \xch $
\xEH
$(\xbm PR)+(\xbm RatM)$
\xEP

\hline

(12.1)
\xEH
$(\xbm =)$
\xEH
$ \xch $ $(\xcs)+(\xbm \xcc)$
\xEH
$(\xbm =')$
\xEP
\cline{1-1}
\cline{3-3}
(12.2)
\xEH
\xEH
$ \xci $
\xEH
\xEP

\hline

(13)
\xEH
$(\xbm \xcc)+(\xbm =)$
\xEH
$ \xch $ $(\xcv)$
\xEH
$(\xbm \xcv)$
\xEP

\hline

(14)
\xEH
$(\xbm \xcc)+(\xbm \xCQ)+(\xbm =)$
\xEH
$ \xch $ $(\xcv)$
\xEH
$(\xbm \xFO),$ $(\xbm \xcv'),$ $(\xbm CUM)$
\xEP

\hline

(15)
\xEH
$(\xbm \xcc)+(\xbm \xFO)$
\xEH
$ \xch $ $(-)$ of $\xdy$
\xEH
$(\xbm =)$
\xEP

\hline

(16)
\xEH
$(\xbm \xFO)+(\xbm \xbe)+(\xbm PR)+$
\xEH
$ \xch $ $(\xcv)$ + sing.
\xEH
$(\xbm =)$
\xEP
\xEH
$(\xbm \xcc)$
\xEH
\xEH
\xEP

\hline

(17)
\xEH
$(\xbm CUM)+(\xbm =)$
\xEH
$ \xch $ $(\xcv)$ + sing.
\xEH
$(\xbm \xbe)$
\xEP

\hline

(18)
\xEH
$(\xbm CUM)+(\xbm =)+(\xbm \xcc)$
\xEH
$ \xch $ $(\xcv)$
\xEH
$(\xbm \xFO)$
\xEP

\hline

(19)
\xEH
$(\xbm PR)+(\xbm CUM)+(\xbm \xFO)$
\xEH
$ \xch $ sufficient,
\xEH
$(\xbm =)$.
\xEP

\xEH
\xEH
e.g., true in $\xdD_{\xdl}$
\xEH
\xEP

\hline

(20)
\xEH
$(\xbm \xcc)+(\xbm PR)+(\xbm =)$
\xEH
$ \xcH $
\xEH
$(\xbm \xFO)$
\xEP

\hline

(21)
\xEH
$(\xbm \xcc)+(\xbm PR)+(\xbm \xFO)$
\xEH
$ \xcH $ (without $(-)$)
\xEH
$(\xbm =)$
\xEP

\hline

(22)
\xEH
$(\xbm \xcc)+(\xbm PR)+(\xbm \xFO)+$
\xEH
$ \xcH $
\xEH
$(\xbm \xbe)$
\xEP

\xEH
$(\xbm =)+(\xbm \xcv)$
\xEH
\xEH
(thus not
\xEP

\xEH
\xEH
\xEH
representable by
\xEP

\xEH
\xEH
\xEH
ranked structures)
\xEP

\hline

\end{tabular}
}
\end{table}

\efa

$ \xCO $

$ \xCO $
\index{Proposition Alg-Log}

\bp

$\hspace{0.01em}$


\label{Proposition Alg-Log}

The following table
``Logical and algebraic rules'' is to be read as follows:

Let a logic $ \xcn $ satisfy $ \xCf (LLE)$ and $ \xCf (CCL),$ and define a
function $f: \xdD_{ \xdl } \xcp \xdD_{ \xdl }$
by $f(M(T)):=M(\ol{ \ol{T} }).$ Then $f$ is well defined, satisfies $(
\xbm dp),$ and $ \ol{ \ol{T} }=Th(f(M(T))).$

If $ \xcn $ satisfies a rule in the left hand side,
then - provided the additional properties noted in the middle for $ \xch $
hold, too -
$f$ will satisfy the property in the right hand side.

Conversely, if $f: \xdy \xcp \xdp (M_{ \xdl })$ is a function, with $
\xdD_{ \xdl } \xcc \xdy,$ and we define a logic
$ \xcn $ by $ \ol{ \ol{T} }:=Th(f(M(T))),$ then $ \xcn $ satisfies $ \xCf
(LLE)$ and $ \xCf (CCL).$
If $f$ satisfies $(\xbm dp),$ then $f(M(T))=M(\ol{ \ol{T} }).$

If $f$ satisfies a property in the right hand side,
then - provided the additional properties noted in the middle for $ \xci $
hold, too -
$ \xcn $ will satisfy the property in the left hand side.

If ``$T= \xbf $'' is noted in the table, this means that, if one of the
theories
(the one named the same way in Definition \ref{Definition Log-Cond}
(page \pageref{Definition Log-Cond}))
is equivalent to a formula, we do not need $(\xbm dp).$

\begin{table}
\caption{Logical and algebraic rules}
\tabcolsep=0.5pt
{\xssc
\begin{tabular}{|c|c|c|c|}

\hline

\multicolumn{4}{|c|}{Basics} \xEP

\hline

(1.1) \xEH $(OR)$ \xEH $\xch$ \xEH $(\xbm OR)$ \xEP

\cline{1-1}

\cline{3-3}

(1.2) \xEH \xEH $\xci$ \xEH \xEP

\hline

(2.1) \xEH $(disjOR)$ \xEH $\xch$ \xEH $(\xbm disjOR)$ \xEP

\cline{1-1}

\cline{3-3}

(2.2) \xEH \xEH $\xci$ \xEH \xEP

\hline

(3.1) \xEH $(wOR)$ \xEH $\xch$ \xEH $(\xbm wOR)$ \xEP

\cline{1-1}

\cline{3-3}

(3.2) \xEH \xEH $\xci$ \xEH \xEP

\hline

(4.1) \xEH $(SC)$ \xEH $\xch$ \xEH $(\xbm \xcc)$ \xEP

\cline{1-1}

\cline{3-3}

(4.2) \xEH \xEH $\xci$ \xEH \xEP

\hline

(5.1) \xEH $(CP)$ \xEH $\xch$ \xEH $(\xbm \xCQ)$ \xEP

\cline{1-1}

\cline{3-3}

(5.2) \xEH \xEH $\xci$ \xEH \xEP

\hline

(6.1) \xEH $(PR)$ \xEH $\xch$ \xEH $(\xbm PR)$ \xEP

\cline{1-1}

\cline{3-3}

(6.2) \xEH \xEH $\xci$ $(\xbm dp)+(\xbm \xcc)$ \xEH \xEP

\cline{1-1}

\cline{3-3}

(6.3) \xEH \xEH $\xcI$ $-(\xbm dp)$ \xEH \xEP

\cline{1-1}

\cline{3-3}

(6.4) \xEH \xEH $\xci$ $(\xbm \xcc)$ \xEH \xEP

\xEH \xEH $T'=\xbf$ \xEH \xEP

\hline

(6.5) \xEH $(PR)$ \xEH $\xci$ \xEH $(\xbm PR')$ \xEP

\xEH \xEH $T'=\xbf$ \xEH \xEP

\hline

(7.1) \xEH $(CUT)$ \xEH $\xch$ \xEH $(\xbm CUT)$ \xEP

\cline{1-1}

\cline{3-3}

(7.2) \xEH \xEH $\xci$ \xEH \xEP

\hline

\multicolumn{4}{|c|}{Cumulativity} \xEP

\hline

(8.1) \xEH $(CM)$ \xEH $\xch$ \xEH $(\xbm CM)$ \xEP

\cline{1-1}

\cline{3-3}

(8.2) \xEH \xEH $\xci$ \xEH \xEP

\hline

(9.1) \xEH $(ResM)$ \xEH $\xch$ \xEH $(\xbm ResM)$ \xEP

\cline{1-1}

\cline{3-3}

(9.2) \xEH \xEH $\xci$ \xEH \xEP

\hline

(10.1) \xEH $(\xcc \xcd)$ \xEH $\xch$ \xEH $(\xbm \xcc \xcd)$ \xEP

\cline{1-1}

\cline{3-3}

(10.2) \xEH \xEH $\xci$ \xEH \xEP

\hline

(11.1) \xEH $(CUM)$ \xEH $\xch$ \xEH $(\xbm CUM)$ \xEP

\cline{1-1}

\cline{3-3}

(11.2) \xEH \xEH $\xci$ \xEH \xEP

\hline

\multicolumn{4}{|c|}{Rationality} \xEP

\hline

(12.1) \xEH $(RatM)$ \xEH $\xch$ \xEH $(\xbm RatM)$ \xEP

\cline{1-1}

\cline{3-3}

(12.2) \xEH \xEH $\xci$ $(\xbm dp)$ \xEH \xEP

\cline{1-1}

\cline{3-3}

(12.3) \xEH \xEH $\xcI$ $-(\xbm dp)$ \xEH \xEP

\cline{1-1}

\cline{3-3}

(12.4) \xEH \xEH $\xci$ \xEH \xEP

\xEH \xEH $T=\xbf$ \xEH \xEP

\hline

(13.1) \xEH $(RatM=)$ \xEH $\xch$ \xEH $(\xbm =)$ \xEP

\cline{1-1}

\cline{3-3}

(13.2) \xEH \xEH $\xci$ $(\xbm dp)$ \xEH \xEP

\cline{1-1}

\cline{3-3}

(13.3) \xEH \xEH $\xcI$ $-(\xbm dp)$ \xEH \xEP

\cline{1-1}

\cline{3-3}

(13.4) \xEH \xEH $\xci$ \xEH \xEP

\xEH \xEH $T=\xbf$ \xEH \xEP

\hline

(14.1) \xEH $(Log =')$ \xEH $\xch$ \xEH $(\xbm =')$ \xEP

\cline{1-1}
\cline{3-3}

(14.2) \xEH \xEH $\xci$ $(\xbm dp)$ \xEH \xEP

\cline{1-1}
\cline{3-3}

(14.3) \xEH \xEH $\xcI$ $-(\xbm dp)$ \xEH \xEP

\cline{1-1}
\cline{3-3}

(14.4) \xEH \xEH $\xci$ $T=\xbf$ \xEH \xEP

\hline

(15.1) \xEH $(Log \xFO)$ \xEH $\xch$ \xEH $(\xbm \xFO)$ \xEP

\cline{1-1}
\cline{3-3}

(15.2) \xEH \xEH $\xci$ \xEH \xEP

\hline

(16.1)
\xEH
$(Log \xcv)$
\xEH
$\xch$ $(\xbm \xcc)+(\xbm =)$
\xEH
$(\xbm \xcv)$
\xEP

\cline{1-1}
\cline{3-3}

(16.2) \xEH \xEH $\xci$ $(\xbm dp)$ \xEH \xEP

\cline{1-1}
\cline{3-3}

(16.3) \xEH \xEH $\xcI$ $-(\xbm dp)$ \xEH \xEP

\hline

(17.1)
\xEH
$(Log \xcv')$
\xEH
$\xch$ $(\xbm \xcc)+(\xbm =)$
\xEH
$(\xbm \xcv')$
\xEP

\cline{1-1}
\cline{3-3}

(17.2) \xEH \xEH $\xci$ $(\xbm dp)$ \xEH \xEP

\cline{1-1}
\cline{3-3}

(17.3) \xEH \xEH $\xcI$ $-(\xbm dp)$ \xEH \xEP

\hline

\end{tabular}
}
\end{table}

\ep

$ \xCO $

$ \xCO $
\index{Definition Preferential structure}

\bd

$\hspace{0.01em}$


\label{Definition Pref-Str}

Fix $U \xEd \xCQ,$ and consider arbitrary $X.$
Note that this $X$ has not necessarily anything to do with $U,$ or $ \xdu
$ below.
Thus, the functions $ \xbm_{ \xdm }$ below are in principle functions from
$V$ to $V$ - where $V$
is the set theoretical universe we work in.

Note that we work here often with copies of elements (or models).
In other areas of logic, most authors work with valuation functions. Both
definitions - copies or valuation functions - are equivalent, a copy
$ \xBc x,i \xBe $ can be seen as a state $ \xBc x,i \xBe $ with valuation $x.$
In the
beginning
of research on preferential structures, the notion of copies was widely
used, whereas e.g., [KLM90] used that of valuation functions. There is
perhaps
a weak justification of the former terminology. In modal logic, even if
two states have the same valid classical formulas, they might still be
distinguishable by their valid modal formulas. But this depends on the
fact
that modality is in the object language. In most work on preferential
stuctures, the consequence relation is outside the object language, so
different states with same valuation are in a stronger sense copies of
each other.

\ed

 \xEh

 \xDH $ \xCf Preferential$ $ \xCf models$ or $ \xCf structures.$
 \index{preferential model}
 \index{preferential structure}

 \xEh

 \xDH The version without copies:

A pair $ \xdm:= \xBc U, \xeb  \xBe $ with $U$ an arbitrary set, and $ \xeb $ an
arbitrary binary relation
on $U$ is called a $ \xCf preferential$ $ \xCf model$ or $ \xCf
structure.$

 \xDH The version with copies \index{copies}:

A pair $ \xdm:= \xBc  \xdu, \xeb  \xBe $ with $ \xdu $ an arbitrary set of
pairs,
and $ \xeb $ an arbitrary binary
relation on $ \xdu $ is called a $ \xCf preferential$ $ \xCf model$ or $
\xCf structure.$

If $ \xBc x,i \xBe  \xbe \xdu,$ then $x$ is intended to be an element of $U,$
and
$i$ the index of the
copy.

We sometimes also need copies of the relation $ \xeb.$ We will then
replace $ \xeb $
by one or several arrows $ \xba $ attacking non-minimal elements, e.g., $x
\xeb y$ will
be written $ \xba:x \xcp y$ \index{$ \xba:x \xcp y$}, $ \xBc x,i \xBe  \xeb
\xBc y,i \xBe $ will
be written
$ \xba: \xBc x,i \xBe  \xcp  \xBc y,i \xBe $ \index{$ \xba: \xBc x,i \xBe  \xcp
\xBc y,i \xBe $}, and
finally we might have $ \xBc  \xba,k \xBe:x \xcp y$ \index{$ \xBc  \xba,k \xBe
:x \xcp y$}  and
$ \xBc  \xba,k \xBe: \xBc x,i \xBe  \xcp  \xBc y,i \xBe $ \index{$ \xBc  \xba,k
\xBe: \xBc x,i \xBe  \xcp  \xBc y,i \xBe $}, etc.

 \xEj

 \xDH $ \xCf Minimal$ $ \xCf elements,$ the functions $ \xbm_{ \xdm }$ \index{$
\xbm_{ \xdm }$}
 \index{minimal element}

 \xEh

 \xDH The version without copies:

Let $ \xdm:= \xBc U, \xeb  \xBe,$ and define

$ \xbm_{ \xdm }(X)$ $:=$ $\{x \xbe X:$ $x \xbe U$ $ \xcu $ $ \xCN \xcE x'
\xbe X \xcs U.x' \xeb x\}.$

$ \xbm_{ \xdm }(X)$ is called the set of $ \xCf minimal$ $ \xCf elements$
of $X$ (in $ \xdm).$

Thus, $ \xbm_{ \xdm }(X)$ is the set of elements such that there is no
smaller one
in $X.$

 \xDH The version with copies:

Let $ \xdm:= \xBc  \xdu, \xeb  \xBe $ be as above. Define

$ \xbm_{ \xdm }(X)$ $:=$ $\{x \xbe X:$ $ \xcE  \xBc x,i \xBe  \xbe \xdu. \xCN
\xcE
\xBc x',i'  \xBe  \xbe \xdu (x' \xbe X$ $ \xcu $ $ \xBc x',i'  \xBe' \xeb  \xBc
x,i \xBe)\}.$

Thus, $ \xbm_{ \xdm }(X)$ is the projection on the first coordinate of the
set of elements
such that there is no smaller one in $X.$

Again, by abuse of language, we say that $ \xbm_{ \xdm }(X)$ is the set of
$ \xCf minimal$ $ \xCf elements$
of $X$ in the structure. If the context is clear, we will also write just
$ \xbm.$

We sometimes say that $ \xBc x,i \xBe $
``$ \xCf kills$'' or ``$ \xCf minimizes$'' $ \xBc y,j \xBe $ if
 \index{kill}
 \index{minimize}
$ \xBc x,i \xBe  \xeb  \xBc y,j \xBe.$ By abuse of language we also say a set
$X$ $ \xCf
kills$ or $ \xCf minimizes$ a set
$Y$ if for all $ \xBc y,j \xBe  \xbe \xdu,$ $y \xbe Y$ there is $ \xBc x,i \xBe
\xbe \xdu,$
$x \xbe X$ s.t. $ \xBc x,i \xBe  \xeb  \xBc y,j \xBe.$

$ \xdm $ is also called $ \xCf injective$ or 1-copy \index{1-copy},
iff there is always at most one
copy
 \index{injective}
$ \xBc x,i \xBe $ for each $x.$ Note that the existence of copies corresponds to
a
non-injective labelling function - as is often used in nonclassical
logic, e.g., modal logic.

 \xEj

 \xEj

We say that $ \xdm $ is $ \xCf transitive,$ $ \xCf irreflexive,$ etc., iff
$ \xeb $ is.
 \index{transitive}
 \index{irreflexive}

Note that $ \xbm (X)$ might well be empty, even if $X$ is not.

$ \xCO $

$ \xCO $
\index{Definition Preferential logics}

\bd

$\hspace{0.01em}$


\label{Definition Pref-Log}

We define the consequence relation \index{consequence relation}  of a
preferential
structure for a
given propositional language $ \xdl.$

 \xEh

 \xDH

 \xEh

 \xDH If $m$ is a classical model of a language $ \xdl,$ we say by abuse
of language

$ \xBc m,i \xBe  \xcm \xbf $ iff $m \xcm \xbf,$

and if $X$ is a set of such pairs, that

$X \xcm \xbf $ iff for all $ \xBc m,i \xBe  \xbe X$ $m \xcm \xbf.$

 \xDH If $ \xdm $ is a preferential structure, and $X$ is a set of $ \xdl
-$models for a
classical propositional language $ \xdl,$ or a set of pairs $ \xBc m,i \xBe,$
where the $m$ are
such models, we call $ \xdm $ a $ \xCf classical$ $ \xCf preferential$ $
\xCf structure$ or $ \xCf model.$

 \xEj

 \xDH

$ \xCf Validity$ in a preferential structure, or the $ \xCf semantical$ $
\xCf consequence$ $ \xCf relation$
defined by such a structure:

Let $ \xdm $ be as above.

We define:

$T \xcm_{ \xdm } \xbf $ \index{$T \xcm_{ \xdm } \xbf $}  iff $ \xbm_{ \xdm
}(M(T)) \xcm
\xbf,$ i.e., $ \xbm_{ \xdm }(M(T)) \xcc M(\xbf).$

 \xDH

$ \xdm $ will be called $ \xCf definability$ $ \xCf preserving$ iff for
all $X \xbe \xdD_{ \xdl }$ $ \xbm_{ \xdm }(X) \xbe \xdD_{ \xdl }.$

 \xEj

As $ \xbm_{ \xdm }$ is defined on $ \xdD_{ \xdl },$ but need by no means
always result in some new
definable set, this is (and reveals itself as a quite strong) additional
property.

\ed

$ \xCO $

$ \xCO $
\index{Definition Smoothness}

\bd

$\hspace{0.01em}$


\label{Definition Smooth}

Let $ \xdy \xcc \xdp (U).$ (In applications to logic, $ \xdy $ will be $
\xdD_{ \xdl }.)$

A preferential structure $ \xdm $ is called $ \xdy -$smooth \index{$ \xdy
-$smooth}  iff for every $X \xbe \xdy $ every
element
$x \xbe X$ is either minimal in $X$ or above an element, which is minimal
in $X.$ More
precisely:

 \xEh

 \xDH The version without copies:

If $x \xbe X \xbe \xdy,$ then either $x \xbe \xbm (X)$ or there is $x'
\xbe \xbm (X).x' \xeb x.$

 \xDH The version with copies:

If $x \xbe X \xbe \xdy,$ and $ \xBc x,i \xBe  \xbe \xdu,$ then either there is
no
$ \xBc x',i'  \xBe  \xbe \xdu,$ $x' \xbe X,$
$ \xBc x',i'  \xBe  \xeb  \xBc x,i \xBe $ or there is $ \xBc x',i'  \xBe  \xbe
\xdu,$
$ \xBc x',i'  \xBe  \xeb  \xBc x,i \xBe,$ $x' \xbe X,$ s.t. there is
no $ \xBc x'',i''  \xBe  \xbe \xdu,$ $x'' \xbe X,$
with $ \xBc x'',i''  \xBe  \xeb  \xBc x',i'  \xBe.$

(Writing down all details here again might make it easier to read
applications
of the definition later on.)

 \xEj

When considering the models of a language $ \xdl,$ $ \xdm $ will be
called $ \xCf smooth$ iff
 \index{smooth}
it is $ \xdD_{ \xdl }-$smooth \index{$ \xdD_{ \xdl }-$smooth} ; $ \xdD_{ \xdl }$
is the
default.

Obviously, the richer the set $ \xdy $ is, the stronger the condition $
\xdy -$smoothness
will be.

\ed

$ \xCO $

$ \xCO $

\bfa

$\hspace{0.01em}$


\label{Fact Rank-Base}

Let $ \xeb $ be an irreflexive, binary relation on $X,$ then the following
two conditions
are equivalent:

(1) There is $ \xbO $ and an irreflexive, total, binary relation $ \xeb'
$ on $ \xbO $ and a
function $f:X \xcp \xbO $ s.t. $x \xeb y$ $ \xcj $ $f(x) \xeb' f(y)$ for
all $x,y \xbe X.$

(2) Let $x,y,z \xbe X$ and $x \xcT y$ wrt. $ \xeb $ (i.e., neither $x \xeb
y$ nor $y \xeb x),$ then $z \xeb x$ $ \xch $ $z \xeb y$
and $x \xeb z$ $ \xch $ $y \xeb z.$

\efa

$ \xCO $

$ \xCO $
\index{Definition Ranked relation}

\bd

$\hspace{0.01em}$


\label{Definition Rank-Rel}

We call an irreflexive, binary relation $ \xeb $ on $X,$ which satisfies
(1)
(equivalently (2)) of Fact \ref{Fact Rank-Base} (page \pageref{Fact Rank-Base})
, ranked \index{ranked}.
By abuse of language, we also call a preferential structure $ \xBc X, \xeb  \xBe
$
ranked, iff
$ \xeb $ is.

\ed

$ \xCO $

$ \xCO $

\bfa

$\hspace{0.01em}$


\label{Fact Rank-Trans}

If $ \xeb $ on $X$ is ranked, and free of cycles, then $ \xeb $ is
transitive.

\efa

$ \xCO $

\subparagraph{
Proof
}

$\hspace{0.01em}$


Let $x \xeb y \xeb z.$ If $x \xcT z,$ then $y \xee z,$ resulting in a
cycle of length 2. If $z \xeb x,$ then
we have a cycle of length 3. So $x \xeb z.$ $ \xcz $
\\[3ex]

$ \xCO $

$ \xCO $

\br

$\hspace{0.01em}$


\label{Remark RatM=}

Note that $(\xbm =')$ is very close to $ \xCf (RatM):$ $ \xCf (RatM)$
says:
$ \xba \xcn \xbb,$ $ \xba \xcN \xCN \xbg $ $ \xch $ $ \xba \xcu \xbg \xcn
\xbb.$ Or, $f(A) \xcc B,$ $f(A) \xcs C \xEd \xCQ $ $ \xch $
$f(A \xcs C) \xcc B$ for all $A,B,C.$ This is not quite, but almost: $f(A
\xcs C) \xcc f(A) \xcs C$
(it depends how many $B$ there are, if $f(A)$ is some such $B,$ the fit is
perfect).

\er

$ \xCO $

$ \xCO $

\bfa

$\hspace{0.01em}$


\label{Fact Rank-Hold}

In all ranked structures, $(\xbm \xcc),$ $(\xbm =),$ $(\xbm PR),$ $(
\xbm ='),$ $(\xbm \xFO),$ $(\xbm \xcv),$ $(\xbm \xcv'),$
$(\xbm \xbe),$ $(\xbm RatM)$ will hold, if the corresponding closure
conditions are
satisfied.

\efa

$ \xCO $

\subparagraph{
Proof
}

$\hspace{0.01em}$


$(\xbm \xcc)$ and $(\xbm PR)$ hold in all preferential structures.

$(\xbm =)$ and $(\xbm =')$ are trivial.

$(\xbm \xcv)$ and $(\xbm \xcv'):$ All minimal copies of elements in
$f(Y)$ have the same rank.
If some $y \xbe f(Y)$ has all its minimal copies killed by an element $x
\xbe X,$ by
rankedness, $x$ kills the rest, too.

$(\xbm \xbe):$ If $f(\{a\})= \xCQ,$ we are done. Take the minimal
copies of $ \xCf a$ in $\{a\},$ they are
all killed by one element in $X.$

$(\xbm \xFO):$ Case $f(X)= \xCQ:$ If below every copy of $y \xbe Y$
there is a copy of some $x \xbe X,$
then $f(X \xcv Y)= \xCQ.$ Otherwise $f(X \xcv Y)=f(Y).$ Suppose now $f(X)
\xEd \xCQ,$ $f(Y) \xEd \xCQ,$ then
the minimal ranks decide: if they are equal, $f(X \xcv Y)=f(X) \xcv f(Y),$
etc.

$(\xbm RatM):$ Let $X \xcc Y,$ $y \xbe X \xcs f(Y) \xEd \xCQ,$ $x \xbe
f(X).$ By rankedness, $y \xeb x,$ or
$y \xcT x,$ $y \xeb x$ is impossible, as $y \xbe X,$ so $y \xcT x,$ and $x
\xbe f(Y).$

$ \xcz $
\\[3ex]

$ \xCO $

$ \xCO $
\index{Proposition Pref-Representation-Without-Ref}

The following table summarizes representation by
preferential structures.

``singletons'' means that the domain must contain all singletons,
``1 copy''
or ``$ \xcg 1$ copy'' means that the structure may contain only 1 copy for
each point,
or several, ``$(\xbm \xCQ)$'' etc. for the preferential structure mean
that the
$ \xbm -$function of the structure has to satisfy this property.

Note that the following table is one (the more difficult) half of a full
representation result for preferential structures. It shows equivalence
between certain abstract conditions for model choice functions and certain
preferential structures. The other half - equivalence between certain
logical rules and certain abstract conditions for model choice functions -
are summarized
in Definition \ref{Definition Log-Cond} (page \pageref{Definition
Log-Cond})  and
shown in Proposition \ref{Proposition Alg-Log} (page \pageref{Proposition
Alg-Log}).

\label{Proposition Pref-Representation-Without-Ref}

\begin{table}[h]

\index{$(\xbm \xcc)$}
\index{reactive}
\index{$(LLE)$}
\index{$(CCL)$}
\index{$(SC)$}
\index{$(\xbm CUM)$}
\index{$(\xbm \xcc \xcd)$}
\index{essentially smooth}
\index{$(SC)$}
\index{$(\xcc \xcd)$}
\index{$(\xbm PR)$}
\index{$(RW)+$}
\index{$(PR)$}
\index{$(\xbm dp)$}
\index{normal characterization}
\index{exception set}
\index{smooth}
\index{$(CUM)$}
\index{$(\xbm=)$}
\index{ranked}
\index{$\xcg 1$ copy}
\index{$(\xbm=')$}
\index{$(\xbm\xFO)$}
\index{$(\xbm\xcv)$}
\index{$(\xbm\xcv')$}
\index{$(\xbm\xbe)$}
\index{$(\xbm RatM)$}
\index{$(\xbm \xCQ)$}
\index{1 copy}
\index{$(\xbm \xCQ fin)$}
\index{$(\xcv)$}
\index{singletons}
\index{$(RatM)$}
\index{$(RatM=)$}
\index{$(Log\xcv)$}
\index{$(Log\xcv')$}

\caption{Preferential representation}

\tabcolsep=0.5pt
\begin{turn}{90}
{\tiny
\begin{tabular}{|c|c|c|c|c|}

\hline

$\xbm-$ function
\xEH
\xEH
Pref.Structure
\xEH
\xEH
Logic
\xEP

\hline

$(\xbm \xcc)+(\xbm PR)$
\xEH
$\xci$
\xEH
general
\xEH
$\xch$ $(\xbm dp)$
\xEH
$(LLE)+(RW)+$
\xEP

\xEH
\xEH
\xEH
\xEH
$(SC)+(PR)$
\xEP

\cline{2-2}
\cline{4-4}

\xEH
$\xch$
\xEH
\xEH
$\xci$
\xEH
\xEP

\cline{2-2}
\cline{4-4}

\xEH
\xEH
\xEH
$\xcH$ without $(\xbm dp)$
\xEH
\xEP

\cline{4-5}

\xEH
\xEH
\xEH
$\xcJ$ without $(\xbm dp)$
\xEH
any ``normal''
\xEP

\xEH
\xEH
\xEH
\xEH
characterization
\xEP

\xEH
\xEH
\xEH
\xEH
of any size
\xEP

\hline

$(\xbm \xcc)+(\xbm PR)$
\xEH
$\xci$
\xEH
transitive
\xEH
$\xch$ $(\xbm dp)$
\xEH
$(LLE)+(RW)+$
\xEP

\xEH
\xEH
\xEH
\xEH
$(SC)+(PR)$
\xEP

\cline{2-2}
\cline{4-4}

\xEH
$\xch$
\xEH
\xEH
$\xci$
\xEH
\xEP

\cline{2-2}
\cline{4-4}

\xEH
\xEH
\xEH
$\xcH$ without $(\xbm dp)$
\xEH
\xEP

\cline{4-5}

\xEH
\xEH
\xEH
$\xcj$ without $(\xbm dp)$
\xEH
using ``small''
\xEP

\xEH
\xEH
\xEH
\xEH
exception sets
\xEP

\hline

$(\xbm \xcc)+(\xbm PR)+(\xbm CUM)$
\xEH
$\xci$
\xEH
smooth
\xEH
$\xch$ $(\xbm dp)$
\xEH
$(LLE)+(RW)+$
\xEP

\xEH
\xEH
\xEH
\xEH
$(SC)+(PR)+$
\xEP

\xEH
\xEH
\xEH
\xEH
$(CUM)$
\xEP

\cline{2-2}
\cline{4-4}

\xEH
$\xch$ $(\xcv)$
\xEH
\xEH
$\xci$ $(\xcv)$
\xEH
\xEP

\cline{2-2}
\cline{4-4}

\xEH
\xEH
\xEH
$\xcH$ without $(\xbm dp)$
\xEH
\xEP

\hline

$(\xbm \xcc)+(\xbm PR)+(\xbm CUM)$
\xEH
$\xci$
\xEH
smooth+transitive
\xEH
$\xch$ $(\xbm dp)$
\xEH
$(LLE)+(RW)+$
\xEP

\xEH
\xEH
\xEH
\xEH
$(SC)+(PR)+$
\xEP

\xEH
\xEH
\xEH
\xEH
$(CUM)$
\xEP

\cline{2-2}
\cline{4-4}

\xEH
$\xch$ $(\xcv)$
\xEH
\xEH
$\xci$ $(\xcv)$
\xEH
\xEP

\cline{2-2}
\cline{4-4}

\xEH
\xEH
\xEH
$\xcH$ without $(\xbm dp)$
\xEH
\xEP

\cline{4-5}

\xEH
\xEH
\xEH
$\xcj$ without $(\xbm dp)$
\xEH
using ``small''
\xEP

\xEH
\xEH
\xEH
\xEH
exception sets
\xEP

\hline

$(\xbm\xcc)+(\xbm=)+(\xbm PR)+$
\xEH
$\xci$
\xEH
ranked, $\xcg 1$ copy
\xEH
\xEH
\xEP

$(\xbm=')+(\xbm\xFO)+(\xbm\xcv)+$
\xEH
\xEH
\xEH
\xEH
\xEP

$(\xbm\xcv')+(\xbm\xbe)+(\xbm RatM)$
\xEH
\xEH
\xEH
\xEH
\xEP

\hline

$(\xbm\xcc)+(\xbm=)+(\xbm PR)+$
\xEH
$\xcH$
\xEH
ranked
\xEH
\xEH
\xEP

$(\xbm\xcv)+(\xbm\xbe)$
\xEH
\xEH
\xEH
\xEH
\xEP

\hline

$(\xbm\xcc)+(\xbm=)+(\xbm \xCQ)$
\xEH
$\xcj$, $(\xcv)$
\xEH
ranked,
\xEH
\xEH
\xEP

\xEH
\xEH
1 copy + $(\xbm \xCQ)$
\xEH
\xEH
\xEP

\hline

$(\xbm\xcc)+(\xbm=)+(\xbm \xCQ)$
\xEH
$\xcj$, $(\xcv)$
\xEH
ranked, smooth,
\xEH
\xEH
\xEP

\xEH
\xEH
1 copy + $(\xbm \xCQ)$
\xEH
\xEH
\xEP

\hline

$(\xbm\xcc)+(\xbm=)+(\xbm \xCQ fin)+$
\xEH
$\xcj$, $(\xcv)$, singletons
\xEH
ranked, smooth,
\xEH
\xEH
\xEP

$(\xbm\xbe)$
\xEH
\xEH
$\xcg$ 1 copy + $(\xbm \xCQ fin)$
\xEH
\xEH
\xEP

\hline

$(\xbm\xcc)+(\xbm PR)+(\xbm \xFO)+$
\xEH
$\xcj$, $(\xcv)$, singletons
\xEH
ranked
\xEH
$\xcH$ without $(\xbm dp)$
\xEH
$(RatM), (RatM=)$,
\xEP

$(\xbm \xcv)+(\xbm\xbe)$
\xEH
\xEH
$\xcg$ 1 copy
\xEH
\xEH
$(Log\xcv), (Log\xcv')$
\xEP

\cline{4-5}

\xEH
\xEH
\xEH
$\xcJ$ without $(\xbm dp)$
\xEH
any ``normal''
\xEP

\xEH
\xEH
\xEH
\xEH
characterization
\xEP

\xEH
\xEH
\xEH
\xEH
of any size
\xEP

\hline

\end{tabular}
}
\end{turn}
\end{table}

$ \xCO $

$ \xCO $
\index{Definition $1-\xca$}

\bd

$\hspace{0.01em}$


\label{Definition 1-infin}

Let $ \xdz = \xBc  \xdx, \xeb  \xBe $ be a preferential structure. Call
$ \xdz $ $1- \xca $ \index{$1- \xca $}  over $Z,$
iff for all $x \xbe Z$ there are exactly one or infinitely many copies of
$x,$ i.e.,
for all $x \xbe Z$ $\{u \xbe \xdx:$ $u= \xBc x,i \xBe $ for some $i\}$ has
cardinality 1 or $ \xcg \xbo.$

\ed

$ \xCO $

$ \xCO $

\bl

$\hspace{0.01em}$


\label{Lemma 1-infin}

Let $ \xdz = \xBc  \xdx, \xeb  \xBe $ be a preferential structure and
$f: \xdy \xcp \xdp (Z)$ with $ \xdy \xcc \xdp (Z)$ be represented by $
\xdz,$ i.e., for $X \xbe \xdy $ $f(X)= \xbm_{ \xdz }(X),$
and $ \xdz $ be ranked and free of cycles. Then there is a structure $
\xdz' $, $1- \xca $ over
$Z,$ ranked and free of cycles, which also represents $f.$

\el

$ \xCO $

\subparagraph{
Proof
}

$\hspace{0.01em}$


We construct $ \xdz' = \xBc  \xdx', \xeb'  \xBe.$

Let $A:=\{x \xbe Z$: there is some $ \xBc x,i \xBe  \xbe \xdx,$ but for all $
\xBc x,i \xBe
\xbe \xdx $ there is
$ \xBc x,j \xBe  \xbe \xdx $ with $ \xBc x,j \xBe  \xeb  \xBc x,i \xBe \},$

let $B:=\{x \xbe Z$: there is some $ \xBc x,i \xBe  \xbe \xdx,$ s.t. for no $
\xBc x,j \xBe
\xbe \xdx $ $ \xBc x,j \xBe  \xeb  \xBc x,i \xBe \},$

let $C:=\{x \xbe Z$: there is no $ \xBc x,i \xBe  \xbe \xdx \}.$

Let $c_{i}:i< \xbk $ be an enumeration of $C.$ We introduce for each such
$c_{i}$ $ \xbo $ many
copies $ \xBc c_{i},n \xBe:n< \xbo $ into $ \xdx',$ put all $ \xBc c_{i},n
\xBe $ above all
elements in $ \xdx,$ and order
the $ \xBc c_{i},n \xBe $ by $ \xBc c_{i},n \xBe  \xeb'  \xBc c_{i' },n'  \xBe
$ $: \xcj $ $(i=i' $ and
$n>n')$ or $i>i'.$ Thus, all $ \xBc c_{i},n \xBe $ are
comparable.

If $a \xbe A,$ then there are infinitely many copies of a in $ \xdx,$ as
$ \xdx $ was
cycle-free, we put them all into $ \xdx'.$
If $b \xbe B,$ we choose exactly one such minimal element $ \xBc b,m \xBe $
(i.e.,
there
is no $ \xBc b,n \xBe  \xeb  \xBc b,m \xBe)$ into $ \xdx',$ and omit all other
elements. (For definiteness, assume in all applications $m=0.)$
For all elements from A and $B,$ we take the restriction of the order $
\xeb $ of $ \xdx.$
This is the new structure $ \xdz'.$

Obviously, adding the $ \xBc c_{i},n \xBe $ does not introduce cycles,
irreflexivity
and
rankedness are preserved. Moreover, any substructure of a cycle-free,
irreflexive,
ranked structure also has these properties, so $ \xdz' $ is $1- \xca $
over $Z,$ ranked and
free of cycles.

We show that $ \xdz $ and $ \xdz' $ are equivalent. Let then $X \xcc Z,$
we have to prove
$ \xbm (X)= \xbm' (X)$ $(\xbm:= \xbm_{ \xdz }$, $ \xbm':= \xbm_{
\xdz' }).$

Let $z \xbe X- \xbm (X).$ If $z \xbe C$ or $z \xbe A,$ then $z \xce \xbm'
(X).$ If $z \xbe B,$
let $ \xBc z,m \xBe $ be the chosen element. As $z \xce \xbm (X),$ there is $x
\xbe
X$ s.t. some $ \xBc x,j \xBe  \xeb  \xBc z,m \xBe.$
$x$ cannot be in $C.$ If $x \xbe A,$ then also $ \xBc x,j \xBe  \xeb'  \xBc z,m
\xBe $. If
$x \xbe B,$ then there is some
$ \xBc x,k \xBe $ also in $ \xdx'.$ $ \xBc x,j \xBe  \xeb  \xBc x,k \xBe $ is
impossible. If $ \xBc x,k \xBe
\xeb  \xBc x,j \xBe,$ then $ \xBc z,m \xBe  \xee  \xBc x,k \xBe $
by transitivity. If $ \xBc x,k \xBe  \xcT  \xBc x,j \xBe $, then also $ \xBc z,m
\xBe  \xee  \xBc x,k \xBe $ by
rankedness. In any
case, $ \xBc z,m \xBe  \xee'  \xBc x,k \xBe,$ and thus $z \xce \xbm' (X).$

Let $z \xbe X- \xbm' (X).$ If $z \xbe C$ or $z \xbe A,$ then $z \xce \xbm
(X).$ Let $z \xbe B,$ and some $ \xBc x,j \xBe  \xeb'  \xBc z,m \xBe.$
$x$ cannot be in $C,$ as they were sorted on top, so $ \xBc x,j \xBe $ exists in
$
\xdx $ too and
$ \xBc x,j \xBe  \xeb  \xBc z,m \xBe.$ But if any other $ \xBc z,i \xBe $ is
also minimal in $ \xdz $
among the $ \xBc z,k \xBe,$
then by rankedness also $ \xBc x,j \xBe  \xeb  \xBc z,i \xBe,$ as $ \xBc z,i
\xBe  \xcT  \xBc z,m \xBe,$ so $z
\xce \xbm (X).$ $ \xcz $
\\[3ex]

$ \xCO $

$ \xCO $

We give a generalized abstract nonsense result, taken
from  \cite{LMS01}, which must be part of the folklore:

\bl

$\hspace{0.01em}$


\label{Lemma Abs-Rel-Ext}


Given a set $X$ and a binary relation $R$ on $X,$ there exists a total
preorder (i.e.,
a total, reflexive, transitive relation) $S$ on $X$ that extends $R$ such
that

$ \xcA x,y \xbe X(xSy,ySx \xch xR^{*}y)$

where $R^{*}$ is the reflexive and transitive closure of $R.$

\el

$ \xCO $

\subparagraph{
Proof
}

$\hspace{0.01em}$


Define $x \xDd y$ iff $xR^{*}y$ and $yR^{*}x.$
The relation $ \xDd $ is an equivalence relation.
Let $[x]$ be the equivalence class of $x$ under $ \xDd.$ Define $[x] \xec
[y]$ iff $xR^{*}y.$
The definition of $ \xec $ does not depend on the representatives $x$ and
$y$ chosen.
The relation $ \xec $ on equivalence classes is a partial order.
Let $ \xck $ be any total order on these equivalence classes that extends
$ \xec.$
Define xSy iff $[x] \xck [y].$
The relation $S$ is total (since $ \xck $ is total) and transitive
(since $ \xck $ is transitive) and is therefore a total preorder.
It extends $R$ by the definition of $ \xec $ and the fact that $ \xck $
extends $ \xec.$
Suppose now xSy and ySx. We have $[x] \xck [y]$ and $[y] \xck [x]$
and therefore $[x]=[y]$ by antisymmetry. Therefore $x \xDd y$ and
$xR^{*}y.$
$ \xcz $
\\[3ex]

$ \xCO $

$ \xCO $

\bp

$\hspace{0.01em}$


\label{Proposition Rank-Rep3}

Let $ \xdy \xcc \xdp (U)$ be closed under finite unions.
Then $(\xbm \xcc),$ $(\xbm \xCQ),$ $(\xbm =)$ characterize ranked
structures for which for all
$X \xbe \xdy $ $X \xEd \xCQ $ $ \xch $ $ \xbm_{<}(X) \xEd \xCQ $ hold,
i.e., $(\xbm \xcc),$ $(\xbm \xCQ),$ $(\xbm =)$ hold in such
structures for $ \xbm_{<},$ and if they hold for some $ \xbm,$ we can
find a ranked relation
$<$ on $U$ s.t. $ \xbm = \xbm_{<}.$ Moreover, the structure can be choosen
$ \xdy -$smooth.

\ep

$ \xCO $

\subparagraph{
Proof
}

$\hspace{0.01em}$


Completeness:

Note that by Fact \ref{Fact Mu-Base} (page \pageref{Fact Mu-Base})  $(3)+(4)$
$(\xbm \xFO),$ $(\xbm \xcv),$ $(\xbm \xcv')$ hold.

Define aRb iff $ \xcE A \xbe \xdy (a \xbe \xbm (A),b \xbe A)$ or $a=b.$
$R$ is reflexive and transitive:
Suppose aRb, bRc, let $a \xbe \xbm (A),$ $b \xbe A,$ $b \xbe \xbm (B),$ $c
\xbe B.$ We show $a \xbe \xbm (A \xcv B).$ By
$(\xbm \xFO)$ $a \xbe \xbm (A \xcv B)$ or $b \xbe \xbm (A \xcv B).$
Suppose $b \xbe \xbm (A \xcv B),$ then $ \xbm (A \xcv B) \xcs A \xEd \xCQ
,$
so by $(\xbm =)$ $ \xbm (A \xcv B) \xcs A= \xbm (A),$ so $a \xbe \xbm (A
\xcv B).$

Moreover, $a \xbe \xbm (A),$ $b \xbe A- \xbm (A)$ $ \xch $ $ \xCN (bRa):$
Suppose there is $B$ s.t. $b \xbe \xbm (B),$
$a \xbe B.$ Then by $(\xbm \xcv)$ $ \xbm (A \xcv B) \xcs B= \xCQ,$ and
by $(\xbm \xcv')$ $ \xbm (A \xcv B)= \xbm (A),$ but
$a \xbe \xbm (A) \xcs B,$ $contradiction.$

Let by Lemma \ref{Lemma Abs-Rel-Ext} (page \pageref{Lemma Abs-Rel-Ext})  $S$ be
a total, transitive, reflexive relation on $U$ which extends
$R$ s.t. $xSy,ySx$ $ \xch $ xRy (recall that $R$ is transitive and
reflexive). Define $a<b$
iff aSb, but not bSa. If $a \xcT b$ (i.e., neither $a<b$ nor $b<a),$ then,
by totality of
$S,$ aSb and bSa. $<$ is ranked: If $c<a \xcT b,$ then by transitivity of
$S$ cSb, but
if bSc, then again by transitivity of $S$ aSc. Similarly for $c>a \xcT b.$

$<$ represents $ \xbm $ and is $ \xdy -$smooth: Let $a \xbe A- \xbm (A).$
By $(\xbm \xCQ),$ $ \xcE b \xbe \xbm (A),$ so bRa,
but (by above argument) not aRb, so bSa, but not aSb, so $b<a,$ so $a \xbe
A- \xbm_{<}(A),$
and, as $b$ will then be $<-$minimal (see the next sentence), $<$ is $
\xdy -$smooth. Let
$a \xbe \xbm (A),$ then for all $a' \xbe A$ aRa', so aSa', so there is
no $a' \xbe A$ $a' <a,$ so $a \xbe \xbm_{<}(A).$

$ \xcz $
\\[3ex]

$ \xCO $
\section{
Acknowledgements
}

We thank A.Herzig, Toulouse, and L.v.d.Torre, Luxembourg, for very helpful
discussions.

\end{document}